\documentclass[12pt,reqno]{amsart}
\usepackage{amsmath,amssymb}
\usepackage[dvips]{graphicx,color}
\usepackage{latexsym,amssymb}
\usepackage{mathrsfs}
\usepackage[abbrev]{amsrefs}

\newtheorem{thm}{Theorem}[section]

\newtheorem{prop}[thm]{Proposition}
\newtheorem{lem}[thm]{Lemma}

 \newenvironment{pf}
    {{\noindent \bf Proof. }}{\hfill $\Box$}

\numberwithin{equation}{section}
\numberwithin{thm}{section}

\begin{document}

\begin{center}\large \bf 
Analyticity and large time behavior 
for the Burgers equation with the critical dissipation 
\end{center}

\footnote[0]
{


E-mail: t-iwabuchi@m.tohoku.ac.jp

}
\vskip5mm

\begin{center}
Tsukasa Iwabuchi 

\vskip2mm

Mathematical Institute, 
Tohoku University\\
Sendai 980-8578 Japan

\end{center}

\vskip5mm

\begin{center}
\begin{minipage}{135mm}
\footnotesize
{\sc Abstract. }
This paper is concerned with the Cauchy problem of 
the Burgers equation with the critical dissipation. 
The well-posedness and analyticity in both of the space and the 
time variables are studied based on the frequency decomposition method. 
The large time behavior is revealed for any large initial data. 
As a result, it is shown that 
any smooth and integrable solution is analytic in space and time 
as long as time is positive and behaves like the Poisson kernel 
as time tends to infinity.  
The corresponding results are also obtined for the quasi-geostrophic equation.

\end{minipage}
\end{center}

\section{Introduction}

\quad 

We consider the Cauchy problem for the Burgers equation 
with the fractional Laplacian: 
\begin{equation}\label{eq:B}
\begin{cases}
\partial_t u + \Lambda u + u \partial _x u = 0, 
& \quad t > 0 , x \in \mathbb R, 
\\
u(0,x) = u_0 (x), 
& \quad x \in \mathbb R , 
\end{cases}
\end{equation}
where $\Lambda := \mathcal F ^{-1}|\xi| \mathcal F $. 
We study the local solvability and the analyticity 
for initial data in the largest space among the Sobolev and Besov spaces 
which are related to the scaling invariant spaces  
and the large time 
behavior of solutions for any large initial data. 
It will be shown  that the similar argument also works for 
the quasi-geostrophic equation.

It is known that the equation \eqref{eq:B} is solvable locally 
in time (see~\cite{KNS-2008,MW-2009,DDL-2009}). 
In fact,  
for any initial data in $H^{\frac{1}{2}} (\mathbb R)$ or  
$\dot B^{\frac{1}{p}}_{p,1} (\mathbb R)$ with $ p < \infty$, 
there eixts a unique local solution.  
Also, the analyticity in space is proved in the papers~\cite{KNS-2008,DDL-2009}, 
while that in time variable has been left open up to now. 
The first purpose of this paper is to construct solutions for initial 
data in the Besov space $B^0_{\infty,\infty} (\mathbb R)$, 
which is analytic in both of the space and the time variables for positive time.

Before recalling the results on the large time behavior, 
let us mention about the global regularity. 
Consider the problem with the fractional Laplacian 
of order $\gamma$. 
\begin{equation}\label{eq:B2}
\partial_t u + \Lambda ^\gamma u 
+ u \partial _x u = 0.
\end{equation}
It was proved that 
the value $\gamma = 1$ 
is the threshold for the occurence of singularity in finite 
time or the global regularity 
(see \cite{ADV-2007,DDL-2009,DGV-2003,KNS-2008}). 
In fact, if $\gamma < 1$, 
it is shown that the gradient of the solution blows up in finite time 
for some continuous initial data. 
On the other hand, if $\gamma \geq 1$, 
such singularity does not appear, so that solution always exists 
globally in time. 
We refer to the results on the global regularizing effects 
in the subcritical case $\gamma > 1$ (see~\cite{DGV-2003}), 
and on the non-uniqueness of weak solutions in the supercritical 
case $\gamma < 1$ (see~\cite{AA-2010}) and 
on the global regularizing effects for the 
$n$-dimensionnal Burgers equation 
in the critical case $\gamma = 1$ (see~\cite{CC-2010}). 
We should also mension about papers concerned with the quasi-geostrophic 
equation~\cite{CaVa-2010,KiNa-2009,KNV-2007}, 
where 
the method of \cite{CaVa-2010} is inspired by Di-Giorgi iterative 
estimates, 
the approach of \cite{KiNa-2009} involves control of H\"older norms 
using appropriate test functions, 
and the proof in \cite{KNV-2007} is based on a nonlocal maximum 
principle and to investigate a certain modulus of continuity of solutions. 
The global regularity is not the aim of this paper, but 
we apply their results to guarantee the global existence.

As for the large time behavior, 
Biler-Karch-Woyczynski~\cite{BKW-1999} 
considered the equation with 
the semigroup generated by 
$ \Lambda ^\gamma - \Delta $ 
$(0 < \gamma < 2 )$ 
to study the asymptotic expansion of solutions 
(see also~\cite{BKW-2000,BKW-2001-2,BKW-2001}). 
For the equation \eqref{eq:B2}, 
Karch-Miao-Xu~\cite{KMX-2008} considered 
the subcritical case $1 < \gamma < 2$   
to study that the large time asymptotics is described  
by the rarefaction waves with some condition in the distance. 
Alibaud-Imbert-Karch~\cite{AIK-2010} investigated that 
the entropy solution converges to the self-similar solution 
for the critical case $\gamma = 1$, 
and the nonlinearity is negligible in the asymptotic expansion 
of solutions for the supercritical case $\gamma < 1$. 
In the previous paper~\cite{Iw-2015}, it was proved that 
the solution behaves like the Poisson kernel if initial 
data is integrable and small in the Besov space 
$\dot B^0_{\infty,1} (\mathbb R)$. 
However, for large initial data, 
the large time behavior of smooth solutions has been left open up to now. 
In our theorem below, 
we show that any smooth and integrable solution 
behaves like the Poisson kernel without smallness 
condition for initial data.

\vskip3mm 

Before stating our results, 
let us recall the definition of Besov spaces. 

\vskip3mm 

\noindent 
{\bf Definition} (Besov spaces). 
Let $\{ \psi \} \cup \{ \phi_j \}_{j \in \mathbb Z}$ be such that 
\begin{gather}\notag 
{\rm supp} \, \widehat{\psi} \subset \{ |\xi| \leq 2 \}, 
\quad 
{\rm supp} \, \widehat{\phi_{j}} \subset \{ 2^{j-1} \leq |\xi| \leq 2^{j+1} \} 
\text{ for any } j \in \mathbb Z, 
\\ \label{deco-inhomo1}
\widehat {\psi} (\xi) + \sum _{j = 1} ^\infty \widehat {\phi_j} (\xi) = 1 
\quad \text{for any } \xi \in \mathbb R^d,
\qquad 
\sum _{ j \in \mathbb Z} \widehat{\phi_j} (\xi) = 1 
\quad \text{for any } \xi \in \mathbb R^d \setminus \{ 0 \}.
\end{gather}
For $s \in \mathbb R$ and $1 \leq p,q \leq \infty$, 
we define the Besov spaces as follows. 

\vskip3mm 

\noindent 
{\rm (i)} $B^s_{p,q} (\mathbb R^d)$ is defined by 
$$\displaystyle 
B^s_{p,q} (\mathbb R^d) 
:= \{ u \in \mathcal S'(\mathbb R^d) \, | \, 
  \| u \|_{B^s_{p,q}} < \infty 
  \},
$$
where 
$$
\| u \|_{B^s_{p,q}} := 
\| \psi * u \|_{L^p} 
+ \big\| \big\{ 2^{sj} \| \phi_j * u \|_{L^p} \big\}_{j \in \mathbb N}
  \big\|_{\ell ^q (\mathbb N)} . 
$$
\noindent 
{\rm (ii)} $\dot B^s_{p,q} (\mathbb R^d)$ is defined by 
$$\displaystyle 
\dot B^s_{p,q} (\mathbb R^d) 
:= \{ u \in \mathcal S'(\mathbb R^d) / \mathcal P(\mathbb R^d) \, | \, 
  \| u \|_{\dot B^s_{p,q}} < \infty 
  \},
$$
where $\mathcal P(\mathbb R^d)$ is the set of all polynomials and 
$$
\| u \|_{\dot B^s_{p,q}} := 
 \big\| \big\{ 2^{sj} \| \phi_j * u \|_{L^p} \big\}_{j \in \mathbb Z}
  \big\|_{\ell ^q (\mathbb Z)} .
$$

\vskip3mm 

We also introduce the standard spaces 
$L^r (0,T ; B^s_{p,q} (\mathbb R^d))$ and the Chemin-Lerner spaces 
$\widetilde L^r (0,T ; B^s_{p,q} (\mathbb R^d))$ which is defined by the set of 
all distributions $u$ such that 
$$
\| u \|_{\widetilde L^r (0,T ; B^s_{p,q}) }
:= \| \psi * u \|_{L^r (0,T ; L^p)} 
+ \big\| \big\{ 2^{sj} \| \phi_j * u \|_{L^r (0,T ; L^p)} \big\}_{j \in \mathbb N}
  \big\|_{\ell ^q (\mathbb N)} 
  < \infty . 
$$

\vskip3mm

We also define the Poisson kernel $P_t (x)$. 
\begin{equation}\notag 
P_t (x) 
:= \mathcal F ^{-1} \big[ e^{- t|\xi|} \big] (x) 
= \frac{\Gamma (\frac{d+1}{2}) \, t^{-d}}
       {\pi ^{\frac{d+1}{2}} (1 + |\frac{x}{t}|^2)^{\frac{d+1}{2}} } 
       \quad \text{for } t > 0, x \in \mathbb R^d,
\end{equation}
where $\Gamma (\cdot)$ is the Gamma function. 
The following is our main theorem. 

\begin{thm}\label{thm:Burgers}
Let $u_0 $ be such that 
\begin{equation}\label{B_initial}
u_0 \in  B^0_{\infty,\infty} (\mathbb R) 
\quad \text{ and } \quad 
\lim_{j \to \infty} \| \phi_j * u_0 \|_{L^\infty} = 0. 
\end{equation}
{\rm (i)}
There exists $T > 0$ and a unique solution $u $ of \eqref{eq:B} such that 
\begin{gather}\notag 
u \in C([0,T],  B^0_{\infty,\infty} (\mathbb R)) 
\cap \widetilde L ^\infty (0,T;  B^0_{\infty,\infty}(\mathbb R))
\cap \widetilde L^1 (0,T ;  B^1_{\infty,\infty} (\mathbb R)), 
\\ \notag 
\lim _{ j \to \infty} \| \phi_j * u(t) \|_{L^\infty} = 0 
\quad \text{for any } t > 0 . 
\end{gather}
Furthermore $u(t,x)$ is real analytic in space and time if $t > 0$. 

\vskip1mm 

\noindent 
{\rm (ii)} 
If $u_0 \in L^1 (\mathbb R) $, 
the solution $u$ in {\rm (i)} satisfies 
that $u(t) \in L^1 (\mathbb R)$ for any $t \geq 0$, 
and that for any $1 \leq p \leq \infty$
\begin{equation}\label{largetime}
\lim _{ t \to \infty} 
t^{1-\frac{1}{p}} 
 \Big\| u(t) - P_t \int_{\mathbb R} u_0 (y) dy \Big\|_{L^p} = 0 .
\end{equation}
Furthermore, for any $\alpha > 0$ there exists $C > 0$ such that 
\begin{equation}\label{decay}
\| |\nabla|^\alpha u(t) \|_{L^p} 
\leq C t^{-(1-\frac{1}{p}) - \alpha} 
\quad \text{for any } t \geq 1. 
\end{equation}

\end{thm}

\noindent 
{\bf Remark. } 
The higher order asymptotic expansion in Theorem~1.2 {\rm (ii)} in \cite{Iw-2015} 
are able to be proved also for any large initial data,  
since we already have the decay estimates 
of all regularity norms as \eqref{decay}. Namely, 
the following assertion 
\begin{equation}\notag 
\begin{split}
\lim_{t \to \infty} 
 t^{1 + d (1 - \frac{1}{p})} 
& \Big\| u(t) -  P_t  M
       + \partial_x P_t \int _{\mathbb R} y u_0 (y) dy 
\\
&
       +\frac{1}{2} (\partial _{x} P_t) 
        \int_0 ^t \int_{\mathbb R} u(\tau , y) ^2 dy d\tau 
       +\frac{1}{2} \int_0^t 
         ( \partial_{x} P_{t-\tau} ) * ( MP_{\tau + 1})^2 d\tau
\\
&
       -\frac{1}{2} (\partial _{x} P_t )
         \int_0^t \int_{\mathbb R} (MP_{\tau + 1} (y))^2 dyd\tau
\Big\|_{L^p} 
= 0 
\end{split}
\end{equation}
is true provided that $u_0 \in L^1 (\mathbb R)$, 
$\int_{\mathbb R} |y u_0 (y)| dy < \infty$  
and \eqref{B_initial} is satisfied, 
where $M = \int_{\mathbb R} u_0 (y) dy $.

\vskip3mm

Let us give remarks on the proof. 
Based on the frequecy decomposition 
method (see \cite{MW-2009,WZ-2011}), 
we develop the local solvability to obtain solutions in more 
general spaces which include non-decaying functions, 
while they~\cite{MW-2009,WZ-2011} considered function spaces 
where the Schwartz class is dense. 
In Proposition~\ref{prop:1} below, 
another frequency localized maximum principle 
is established, 
which enables us to work with the iterative method. 
On the analyticity, the existing method is to consider 
the direct computation in the frequency space through 
the Plancherel theorem 
(see e.g.~\cite{DDL-2009,KNS-2008}), 
or to apply the Fourier multiplier theorem 
to multipliers $e^{t\Lambda}$, 
$e^{t (|\partial_{x_1}| + \cdots + |\partial_{x_d}|)}$ 
(see~\cite{BBT-2012,DDL-2009,Lemarie_2002}), 
which requires the boundedness of the Riesz transformation 
or singular integral operators. 
On the other hand, we simply consider 
the derivatives with the weight of time variable, 
which does not require the Fourier multiplier theorem, 
and the analyticity of Gevery type with order one  
is verified for not only space 
but also time variable 
(see Propositions~\ref{prop:2}, \ref{prop:anal_nonlinear} 
and the proof of analyticity thereunder).

As for the large time behavior \eqref{largetime} 
without smallness of initial data, 
the most important point is to get an integrability 
of solutions 
\begin{equation}\label{0805-1}
\int_0^\infty \| u(t) \|_{\dot B^1_{\infty,1}} \, dt < \infty,
\end{equation}
which can be seen in the previous paper~\cite{Iw-2015} and 
is useful to work via the Duhamel formula. 
To handle large initial data, we develop the time decay estimate in 
$L^\infty (\mathbb R)$ along the paper~\cite{CoCo-2004}, 
which assures the smallness of $u(t_0)$ for some large $t_0$. 
Then we can apply the argument for small data for $t \geq t_0$, 
while the integrability for $[0,t_0]$ is guaranteed by 
$u \in L^1_{loc} ([0,\infty) , \dot B^1_{\infty,1} (\mathbb R))$ 
by the result~\cite{MW-2009}. 
In addition, we also establish decay estimates with arbitrary positive regularity 
\eqref{decay} for any large initial data (see also Proposition~\ref{prop:decay}). 

\vskip3mm

We next consider the quasi-geostrophic equation. 
\begin{equation}\label{eq:QG}
\begin{cases}\displaystyle 
 \partial_t \theta + \Lambda ^\gamma \theta  
  + (u \cdot \nabla ) \theta =0,
 & \quad  t > 0 , x \in \mathbb R ^2, 
\\
 u= (-R_2 \theta , R_1 \theta) , 
 & \quad  t > 0 , x \in \mathbb R ^2, 
\\ 
 \theta(0,x) = \theta_0(x) , 
 & \quad x \in \mathbb R ^2, 
\end{cases}
\end{equation}
where $\theta : \mathbb R^2 \to \mathbb R$ 
is the unknown function, 
$R_1, R_2$ are the Riesz transforms in $\mathbb R^2$. 
The equation is an important model in geophysical fluid dynamics, 
which describes the evolution of a surface temperature field 
in a rotating and stratified fluid.  
We obtain the following result by modifiying the method for the Burgers equation 
to handle the Riesz transform especially for the 
low frequency part. For the sake of simplicity, 
we consider initial data in 
$L^1 (\mathbb R^2) \cap B^0_{\infty,\infty} (\mathbb R^2)$ 
to avoid the complexity to get the following result. 

\begin{thm}\label{thm:QG}
Let $\gamma = 1$, $u_0 $ satisfy 
$u_0 \in L^1  (\mathbb R^2) \cap B^0_{\infty, \infty} (\mathbb R^2)$, 
$\| \phi_j * u_0 \|_{L^\infty} \to 0$ as $j \to \infty$. 
Then, there exists 
a unique global solution $u $ of \eqref{eq:QG} such that 
\begin{gather}\notag 
u \in C([0,\infty),  B^0_{\infty,\infty} (\mathbb R^2)) 
\cap \widetilde L ^\infty_{loc} (0,\infty;  B^0_{\infty,\infty}(\mathbb R^2))
\cap \widetilde L^1_{ loc} ([0,\infty) ;  B^1_{\infty,\infty} (\mathbb R^2)), 
\\ \notag 
\lim _{ j \to \infty} \| \phi_j * u(t) \|_{L^\infty} = 0 
\quad \text{for any } t > 0 , 
\end{gather}
and $u(t,x)$ is real analytic in space and time if $t > 0$. 
Furthermore, 
for any $1 \leq p \leq \infty$ 
\begin{equation}\notag 
\lim _{ t \to \infty} 
t^{2(1-\frac{1}{p})} 
 \Big\| u(t) - P_t \int_{\mathbb R^2} u_0 (y) dy \Big\|_{L^p} = 0 , 
\end{equation}
and for any $\alpha > 0$ there exists $C > 0$ such that 
\begin{equation}\notag 
\| |\nabla|^\alpha u(t) \|_{L^p} 
\leq C t^{-2(1-\frac{1}{p}) - \alpha} 
\quad \text{for any } t \geq 1. 
\end{equation}
\end{thm}

\vskip3mm

%

Let us compare with known results. There are a lot of known results 
on the global solvability and the asymptotic behavior. 
The solvability in the case when $\gamma = 1$ was studied for small 
initial data in $L^\infty (\mathbb R^2)$ 
in~\cite{CoCoWu-2001} and that for 
arbitrarily large data has been settled 
in the papers~\cite{CaVa-2010,KiNa-2009,KNV-2007}. 
The well-posedness is also studied in Besov spaces 
in~\cite{WZ-2011}, where spaces are defined by the completion of 
the Schwartz class. 
Here we refer the recent papers~\cite{BMS-2015,CotVic-2016,HmKe-2007,Kise-2011} 
on regularity of super-critical case $\gamma < 1$, where the global regularity is open.
As for the large time behavior,  the subcritical case $\gamma > 1$ 
was resolved by \cite{ConWu-1999,DoLi-2008,SchSch-2003}. 
For the critical case $\gamma = 1$, 
the time decay estimates in $L^p (\mathbb R^2)$ for some $p$ of weak solutions 
are known in \cite{CoCoWu-2001,DoDu-2008,NiSch-2015}. 
We also refer to the recent papers \cite{FeNiPl-2017} for modified equations, 
including not only critical case but also super-critical case, 
\cite{CoTaVi-2015} on the existence of a compact 
global attracter with a time-independent force, and see also references therein. 
The contribution of Theorem \ref{thm:QG} is 
the analyticity in the space and the time variables and to reveal the large time 
behavior of smooth solutions for any large initial data.

\vskip3mm

This paper is organized as follows. 
In section 2, we prepare lemmas on the frequency localized maximum 
principle. 
On the proof of theorems, we prove Theorem~\ref{thm:Burgers} only, 
since Theorem~\ref{thm:QG} follows analogously. 
Section 3 is devoted to showing the local solvability and the analyticity 
for Burgers equation.  
In section 4, we verify the large time behavior of solutions 
in Theorem~\ref{thm:Burgers}.

\section{Preliminary}

In this section, we prepare the frequency localized maximum principle 
for non-decaying smooth functions, which is motivated by the paper~\cite{WZ-2011}. 
The Fourier multiplier theorem for the Poisson kernel and 
the continuity property of linear solutions are also investigated. 

\begin{prop}\label{prop:1} 
{\rm (}Frequency localized maximum principle{\rm)}
Let $u, v, f$ be smooth functions on $(0,\infty) \times \mathbb R^d$ 
such that $u, \partial_t u, \partial_t^2 u, v \in L^\infty (\mathbb R^d)$. 

\noindent 
{\rm (i)} Let $j \in \mathbb Z$. If $u$ satisfies 
$\partial _t \big( \phi_j * u \big) 
+ (v \cdot \nabla ) \big( \phi_j * u \big)
+ \Lambda \big( \phi_j * u \big)= f$, 
then there exists a positive constant $c$ independent of 
$u,v,f, j$ such that for almost every $t \geq 0$ 
\begin{equation}\label{0606-1}
\partial_t \| \phi_j * u \|_{L^\infty} 
+ c 2^j \| \phi_j * u \|_{L^\infty} 
\leq \| f \|_{L^\infty} .
\end{equation}

\noindent
{\rm (ii)} 
If $u$ satisfies 
$\partial _t \big( \psi *u \big) 
+ (v \cdot \nabla ) \big( \psi * u \big)
+ \Lambda \big( \psi * u \big)= f$, 
then for almost every $t \geq 0$ 
\begin{equation}\label{0606-1_inho}
\partial_t \| \psi * u \|_{L^\infty} 
\leq \| f \|_{L^\infty} .
\end{equation}

\end{prop}

In order to prove the above proposition, 
we prepare two lemmas.

\begin{lem}\label{lem:1}
Let $u(t,x)$ be a smooth function  
with $u, \partial_t u, \partial_t^2 u \in L^\infty (\mathbb R^d)$. 
Then $\partial_t \| u(t) \|_{L^\infty}$ exists for almost 
every $t \geq 0$. Furthermore, for almost every $t \geq 0$, 
there exists a sequence $\{ x_{t,n} \}_{n=1}^\infty \subset 
\mathbb R^d$ such that 
\begin{gather}
 \label{0512-2-2}
\|u(t)\|_{L^\infty} = 
\lim_{n \to \infty} 
u (t,x_{t,n}) \operatorname{sgn} \big( u(t,x_{t,n}) \big), \\ 
\label{0512-2}
\partial_t \| u(t) \|_{L^\infty} 
= \lim_{n \to \infty} 
 (\partial_t u) (t,x_{t,n}) \operatorname{sgn} \big( u(t,x_{t,n}) \big), 
\end{gather}
where $\operatorname{sgn} u$ is a sign function of $u$. 
\end{lem}

\begin{pf}
We prove based on the proof of Lemma~3.2 in the paper \cite{WZ-2011}, 
but there needs some modification to handle non-decaying functions. 
In what follows, let $t \geq 0$ and $|h| < 1$ such that $t + h \geq 0$. 

The existence of $\partial_t \| u(t) \|_{L^\infty}$ for almost 
every $t$ is proved by 
\begin{equation}\notag 
\big| \| u(t+h) \|_{L^\infty} - \| u(t) \|_{L^\infty } \big|
\leq 
\| u(t+h) - u(t) \|_{L^\infty}
\leq 
 \sup _{ |\tau -t|\leq 1}\| \partial_\tau u (\tau) \|_{L^\infty} |h|, 
\end{equation}
since this inequality impies the Lipschitz continutity of 
$\| u(t) \|_{L^\infty}$. 

We turn to prove the latter assertion. 
Let us consider non-negative $u$ for the sake of simplicity. 
For each $t+h$, there exists a sequence $\{ x_{t+h,n} \}_{n=1}^\infty$ 
such that 
\begin{equation}\label{0512-9}
u(t+h,x_{t+h,n}) \to \| u(t+h) \|_{L^\infty} 
\quad \text{as } n \to \infty .
\end{equation}
By considering the limit as $h \to 0$, we can take sequences 
$\{ h_m \}_{m=1}^\infty$ and $\{ x_{t+h_m,n_m} \}_{m =1}^\infty$ 
such that 
\begin{gather}\notag 
h_m \to 0 , \quad 
u(t+h_m , x_{t+h_m , n_m}) 
\to \| u(t) \|_{L^\infty} 
\quad \text{as } m \to \infty , 
\\ \label{0512-3}
\| u(t+h_m) \|_{L^\infty} 
- u(t+h_m , x_{t+h_m,n_m} ) 
\leq h_m ^2 ,
\\ \label{0512-3-(2)}
\| u(t) \|_{L^\infty} - u(t , x_{t+h_m, n_m}) \leq h_{n_m} ^2, 
\\ \notag 
\lim_{ m \to \infty} (\partial_t u) (t, x_{t+h_{m},m})  
\text{ exists}.  
\end{gather}
Put 
$
x'_m := x_{t+h_m,n_m}
$.
It follows from the inequality \eqref{0512-3} 
and smoothness of $u$ in the time variable that 
\begin{equation}\notag 
\begin{split}
\frac{\| u(t+h_{m}) \|_{L^\infty} - \| u(t) \|_{L^\infty}}{h_m} 
\leq 
& 
\frac{ u(t+h_{m}, x'_{m} ) + h_{m}^2- u(t, x_{m}')}{h_m} 
\\
\leq 
& 
\frac{ u(t+h_{m}, x'_{m} ) - u(t, x_{m}')}{h_{m}} + h_{m}
\\
\to 
& \lim_{m \to \infty} ( \partial_t u  )(t,x'_m), 
\end{split}
\end{equation}
as $m \to \infty$, which implies 
\begin{equation}\label{0512-6}
\partial_t \| u(t) \|_{L^\infty} \leq 
\lim_{m \to \infty} ( \partial_t u  )(t,x'_m) .
\end{equation}
On the other hand, we have from \eqref{0512-3-(2)} that 
\begin{equation}\notag 
\begin{split}
\frac{\| u(t+h_{m}) \|_{L^\infty} - \| u(t) \|_{L^\infty}}{h_m} 
\geq 
& 
\frac{ u(t+h_{m}, x'_{m} ) - u(t, x_{m}') - h_{m}^2 }{h_m} 
\\
\geq 
& 
\frac{ u(t+h_{m}, x'_{m} ) - u(t, x_{m}')}{h_{m}} - h_{m}
\\
\to 
& \lim_{m \to \infty} ( \partial_t u  )(t,x'_m), 
\end{split}
\end{equation}
as $m \to \infty$, which proves the inequality 
in the opposite direction of \eqref{0512-6}. Hence, 
the assertion \eqref{0512-2} of (a) together with \eqref{0512-2-2} 
is proved for non-negative $u$. 
For general functions $u$, the analogous argument also works well 
by replacing $u$ with $-u$ at which $u$ takes the negative value. 
Therefore we complete the proof. 
\end{pf}

\begin{lem}\label{lem:2}
Let 
$$
\mathcal A := 
\big\{ g \in L^\infty (\mathbb R^d) \,\,\, \big| \, \,\,
\| g \|_{L^\infty} = 1, \,\,\,
{\rm supp \, } \widehat g \subset 
\{ \xi \in  \mathbb R^d  \, | \, 2^{-1} \leq |\xi| \leq 2\} 
\big\}. 
$$
Suppose that $g \in \mathcal A$ and 
$\{ x_n \}_{n=1}^\infty$ satisfies 
$$
\lim_{n \to \infty} g(x_n) \operatorname{sgn} (g(x_n)) 
= \| g \|_{L^\infty} .
$$
Then there exists a positive constant $c$ independent of $g$ 
such that 
\begin{equation}\label{0512-10}
\lim_{n \to \infty} \big(\Lambda g(x_n) \big) \operatorname{sgn} (g(x_n)) 
\geq c.  
\end{equation}
\end{lem}

\begin{pf}
Assume that there exist 
$\{ g_n \}_{n=1} ^\infty$ and $\{ x_n \}_{n=1}^\infty$ 
such that 
\begin{equation}\label{0512-11}
g_n (x_n) \geq 1 - \frac{1}{n}, 
\quad  ( \Lambda g_n ) (x_n) \to 0 \text{ as } n \to \infty ,
\end{equation}
noting that it suffices to consider $-g_n$ instead of $g_n$ 
if $g_n$ can be negative. 
Let 
$$
\tilde g_n (x) := g_n(x+ x_n) . 
$$
Noting that $\tilde g_n(0) \to 1$ as $n \to \infty$ and 
$\| \nabla g_n \|_{L^\infty} \leq C \| g_n \|_{L^\infty} = C$, by taking a 
subsequence we see that 
there exists $g \in L^\infty (\mathbb R^d)$ with 
$g(0) = 1$ such that 
$\tilde g_n(x) \to g(x)$ as $n \to \infty$ uniformly with respect to $x \in \mathbb R^d$ 
by the Ascoli-Arzel\'a theorem 
and $g \in C^\infty (\mathbb R^d) \cap L^\infty (\mathbb R^d)$, 
since $\widehat g_n$ is supported in the bounded set 
$\{ \xi \in \mathbb R^d \, | \,  2^{-1} \leq |\xi| \leq 2 \}$ and so is $g$. 
This $g$ also satisfies $g \not \equiv 1$, 
since 
$
1 = g(0) = \big((\phi_{-1} + \phi_0 + \phi_1) * g \big)(0)
$
but the constant function $ 1$  satisfies  
$(\phi_{-1} + \phi_0 + \phi_1) *1 = 0$.
Here we recall the following formula:
\begin{equation}\label{0608-1}
\Lambda g(x) 
= C_{d} 
\int_{\mathbb R^d} \frac{2g(x)-g(x+y)-g(x-y)}{|x-y|^{d+1}} \,dy, 
\end{equation}
where $C_d$ is a positive constant depending on the dimensions 
(see e.g. Lemma~3.2 in \cite{DPV-2012}).  
The above formula yields that 
\begin{equation}\notag 
\operatorname{sgn} \big( g(0) \big) \Lambda g(0) 
= 
\operatorname{sgn} \big( g(0) \big) C_{d} 
\int_{\mathbb R^d} \frac{2g(0)-g(y)-g(-y)}{|y|^{d+1}} \,dy
> 0 ,
\end{equation}
since $g$ is not a constant function. 
Here we can show that 
\[
\lim_{n \to \infty} 
\int_{\mathbb R^d} 
 \frac{2\tilde g_n(0)-\tilde g_n(y)-\tilde g_n(-y)}{|y|^{d+1}} 
 \,dy
= 
\int_{\mathbb R^d} 
 \frac{2\tilde g(0)-\tilde g(y)-\tilde g(-y)}{|y|^{d+1}} 
 \,dy .
\]
In fact, for any $\delta > 0$, it follows from the dominated convergence theorem that 
\[
\lim_{n \to \infty} 
\int_{ |y| > \delta} 
 \frac{2\tilde g_n(0)-\tilde g_n(y)-\tilde g_n(-y)}{|y|^{d+1}} 
 \,dy
= 
\int_{ |y| > \delta} 
 \frac{2\tilde g(0)-\tilde g(y)-\tilde g(-y)}{|y|^{d+1}} 
 \,dy ,
\]
while the uniform boundedness of 
$\|\nabla ^2 g_n\| _{L^\infty} $ with respect to $n$ 
gives that 
\[
\begin{split}
\Big| 
\int_{ |y| < \delta} 
 \frac{2\tilde g_n(0)-\tilde g_n(y)-\tilde g_n(-y)}{|y|^{d+1}} 
 \,dy
\Big| 
\leq 
& 
\| \nabla ^2 g_n \|_{L^\infty} 
\int _{|y| < \delta} \frac{|y|^2}{|y|^{d+1}} dy 
\\
\leq 
& 
C \| g_n \|_{L^\infty } \delta 
= C \delta . 
\end{split}
\]
Therefore, by the assumption \eqref{0512-11}, we find that 
\begin{equation}\notag 
\begin{split}
0< 
\operatorname{sgn} \big( g(0) \big) \Lambda g(0) 
\leq  
& 
\operatorname{sgn} \big( g(0) \big)
\lim_{n \to \infty} C_{d} 
\int_{\mathbb R^d} 
 \frac{2\tilde g_n(0)-\tilde g_n(y)-\tilde g_n(-y)}{|y|^{d+1}} 
 \,dy
\\
= 
& 
\operatorname{sgn} \big( g(0) \big)
\liminf _{ n \to \infty}\Lambda g_n(x_n) 
\\
= 
& 0 . 
\end{split}
\end{equation}
which is contradiction. 
\end{pf}

\vskip3mm 

We are ready to prove Proposition \ref{prop:1}. 

\vskip3mm 

\noindent 
{\bf Proof of Proposition \ref{prop:1}.}
The inequality \eqref{0606-1} is an immidiate consequence of 
Lemmas~\ref{lem:1}, \ref{lem:2} with the scaling by $2^j$ 
and the fact that 
for any $\{ x_{t,n} \}_{n = 1} ^\infty$ satisfying 
$$
\|u(t)\|_{L^\infty} = 
\lim_{n \to \infty} 
u (t,x_{t,n}) \operatorname{sgn} \big( u(t,x_{t,n}) \big), 
$$
its gradient must satisfy 
$$
\lim_{n \to \infty} \nabla u (t,x_{t,n}) = 0 .
$$
It is also readily to show \eqref{0606-1_inho}, 
since the same argument works well for $\partial _t (\psi * u)$ 
by the non-negativity $0 \leq \liminf_{n \to \infty} 
\Lambda \psi* u(t, x_{t,n}) \operatorname{sgn} (\psi *u(t, x_{t,n})) $ 
for $\{ x_{t,n} \}_{n=1}^\infty$ satisfying 
$|\psi * u(t, x_{t,n})| \to \| \psi * u(t) \|_{L^\infty}$ as $n \to \infty$, 
which is seen from \eqref{0608-1}. 
\hfill $\Box$

\vskip3mm

The following is the Fourier multiplier theorem 
for the propagator defined with the Poisson kernel.

\begin{lem}{\rm (see e.g. \cite{MW-2009,Iw-2015})}
There exist $C>0$, $0 < c < 1$ independent of $u_0$ such that 
\begin{equation}\label{0512-13}
c e^{-Ct2^j} \| \phi_j * u_0 \|_{L^\infty}
\leq 
\| \phi_j * (e^{-t\Lambda} u_0) \|_{L^\infty} 
\leq C e^{-ct2^j} \| \phi_j * u_0 \|_{L^\infty}
\end{equation}
for all $j \in \mathbb Z$ and $u_0 \in L^\infty (\mathbb R^d)$. 

\end{lem}

Next lemma is concerned with the continuity of the linear solution. 

\begin{lem}
{\rm (i)} 
Let $u_0 \in  B^0_{\infty,\infty} (\mathbb R^d)$. Then 
$$
\lim _{t \to 0} e^{-t\Lambda} u_0 = u_0 
\text{ in }   B^0_{\infty,\infty }(\mathbb R^d)
\quad \text{if and only if} \quad 
\lim_{ j \to \infty} \| \phi_j * u_0 \|_{L^\infty} = 0.
$$

\noindent 
{\rm (ii)} 
Let $u_0 \in  B^0_{\infty,\infty} (\mathbb R^d)$ be such that 
$\| \phi_j * u_0 \|_{L^\infty} \to 0 $ as $j \to \infty$. 
Then 
\begin{equation}\label{0512-12}
\lim _{ T \to 0} 
\| e^{-t\Lambda} u_0 \|_{\widetilde L^1 (0,T;  B^1_{\infty,\infty})}
= 0.
\end{equation}

\end{lem}
\begin{pf} 
We prove {\rm (i)} first. Since 
for any $\xi \in \mathbb R^d$ with $2^{j-1} \leq |\xi| \leq 2^{j+1}$ 
$$
| e^{-t|\xi|} -1 | \simeq 
\begin{cases} 
1 & \quad \text{if } t 2^{j} \geq 1,
\\
t2^j & \quad \text{if } t 2^{j} \leq 1,
\end{cases}
$$
we have from the Fourier multiplier theorem that 
\begin{equation}\notag 
\begin{split}
C^{-1} \min \{1, t2^{j}\} \| \phi_j * u_0 \|_{L^\infty}
\leq 
&
\| \phi_j * \big((e^{-t\Lambda} -1 )u_0 \big) \|_{L^\infty} 
\\
\leq 
& C \min\{ 1, t 2^{j}\} \| \phi_j * u_0 \|_{L^\infty} .
\end{split}
\end{equation}
The above first inequality implies the high frequency part of $u_0$ 
vanishing by the time continuity of $e^{-t\Lambda} u_0$, 
and the second one yields 
the time continuity under the condition that high frequency of $u_0$ vanishes. 

Let us turn to prove {\rm (ii)}. 
 For $j_0 \in \mathbb N$, 
it follows from \eqref{0512-13} that 
\begin{equation}\notag 
\begin{split}
\| e^{-t\Lambda} u_0 \|_{\widetilde L^1 (0,T;  B^1_{\infty,\infty})}
\leq 
& T \Big( \| \psi * u_0 \|_{L^\infty} + 
    \sup_{ 1 \leq j \leq j_0 } 2^j \| \phi_j * u_0 \|_{L^\infty}
    \big) 
\\
& +C \sup_{j > j_0} \| e^{-c t 2^{j}} \|_{L^1 (0,T)} 2^j 
    \| \phi_j * u_0 \|_{L^\infty}
\\
\leq 
& T 2^{j_0} \| u_0 \|_{ B^0_{\infty,\infty}} + 
C \sup_{j > j_0} \| \phi_j * u_0 \|_{L^\infty} .
\end{split}
\end{equation}
By taking $j_0$ large and choosing $T$ sufficiently small 
$T \ll 2^{-j_0} \| u_0 \|_{ B^0_{\infty,\infty}}$, 
we get \eqref{0512-12}. 
\end{pf}

\section{Local solvability and analyticity for Burgers equation} 
\label{sec:LWP}

In this section, we prove the local in time solvability and analyticity 
in Theorem~\ref{thm:Burgers}. 
Only in this section, let $\{ \phi_j \}_{j = 0} ^\infty$ be 
the Littlewood Paley dyadic decomposition 
for inhomogeneous spaces, namely, 
$\phi_0$ is taken as $\phi_0 = \psi$, where $\psi$ satisfy 
\eqref{deco-inhomo1}. 
Put 
$S_j := \sum _{k = 0} ^j \phi_j * $ for $j = 0,1,2,\cdots$ 
and $S_j = 0$ for $j = -1,-2,\cdots$ for the sake of simplicity. 
Consider a sequence $\{ u_n \}_{n=1}^\infty$ such that 
\begin{equation}\label{0402-2}
\begin{cases}
u_1 = e^{-t\Lambda} S_1u_0, 
\\
\displaystyle 
\partial _t u_{n+1} + \Lambda u_{n+1} 
+ \sum _{ l \geq 0} 
    \big( S_{ l-3} u _{n} \big) 
    \partial _x \phi_l * u_{n+1} 
\\
\displaystyle 
\quad 
= -  \sum _{ l \geq 0} 
    \Big( \sum_{k \geq l+3} \phi_k * u _{n} 
    \Big) 
    \partial _x \phi_l * u_{n} 
  - \frac{1}{2} \partial _x \sum _{ |l-k| \leq 2} 
     (\phi_k * u_n) (\phi_l * u_n) , 
\\
u_n(0) = S_n u_0.      
& 
\end{cases}
\end{equation}
Existence of $u_{n+1}$ for given $u_n$ is assured by smoothness of 
the initial data $S_n u_0$, and we need to obtain a priori estimate. 
It follows from the boundedness of  
$e^{-t\Lambda}$ and the maximal regularity estimate 
in $ B^0_{\infty,\infty} (\mathbb R)$ that 
$$
\| u_1 \|_{\widetilde L^\infty([0,T] ;  B^0_{\infty,\infty})
        \cap \widetilde L^1 (0,T ;  B^1_{\infty,\infty})}  
+ 
\| \partial _t u_1 \|_{\widetilde L^1 (0,T;  B^0_{\infty,\infty})} 
\leq C \| u_0 \|_{ B^0_{\infty,\infty}} . 
$$
We need to estimate 
$u_n$ $(n = 2,3, \cdots)$ and the difference 
$u_{m+1} - u_{n+1}$. For this purpose, we prepare the 
following lemma.

\begin{lem}\label{lem:3}
Let $T > 0$, $1 \leq q, r \leq \infty$, $s > - 1$ and $0 < \delta < 1/2$. 
Then there exist positive constant $C, c $
such that the following three inequalities 
hold. 
\begin{gather} \label{0402-3}
\begin{split}
\| u_{n+1} \|_{\widetilde{L} ^\infty(0,T; B^{s}_{\infty,q})} 
\leq  
&
\| u_0 \|_{ B^{s}_{\infty,q}} 
+  C  \| u_n \|_{\widetilde L^2 (0,T ;  B^\frac{1}{2}_{\infty,\infty})}
      \| u_{n+1} \|_{\widetilde L^2 (0,T ;  B^{\frac{1}{2}+s}_{\infty,q})}
\\
&+  C  \| u_n \|_{\widetilde L^2 (0,T ;  B^\frac{1}{2}_{\infty,\infty})}
      \| u_n \|_{\widetilde L^2 (0,T ;  B^{\frac{1}{2}+s}_{\infty,q})}, 
\end{split}
\\ 
\label{0402-4}
\begin{split}
\| u_{n+1} \|_{\widetilde L^r (0,T ;  B^{\frac{1}{r}+s}_{\infty,q})} 
\leq 
&
\| e^{-tc \Lambda} u_0 \|_{\widetilde L^r (0,T ;  B^{\frac{1}{r}+s}_{\infty,q})} 
+  C  \| u_n \|_{\widetilde L^2 (0,T ;  B^\frac{1}{2}_{\infty,\infty})}
      \| u_{n+1} \|_{\widetilde L^2 (0,T ;  B^{\frac{1}{2}+s}_{\infty,q})}
\\
& 
+  C  \| u_n \|_{\widetilde L^2 (0,T ;  B^\frac{1}{2}_{\infty,\infty})}
      \| u_n \|_{\widetilde L^2 (0,T ;  B^{\frac{1}{2}+s}_{\infty,q})} , 
\end{split}
\\ 
\label{0402-d}
\begin{split}
\| u_{n+1} - u_{m+1} \|
_{\widetilde L^2 (0,T ;  B^{-\delta}_{\infty,q})}
\leq 
& 
C \| e^{-tc \Lambda} (S_{n+1}u_0 - S_{m+1} u_0) 
  \|_{\widetilde L^2 (0,T ;  B^{-\delta}_{\infty,q})} 
\\
&
+ C 
 \| u_n \|_{\widetilde L^2 (0,T ;  B^{\frac{1}{2}}_{\infty,\infty})}
  \| u_{n+1} - u_{m+1} 
  \|_{\widetilde L^2 (0,T ;  B^{-\delta}_{\infty,q} )}
\\
& 
+ C\big( \| u_{n} \|_{\widetilde L^2 (0,T ;  B^{\frac{1}{2}}_{\infty,\infty})}
    + \| u_{m} \|_{\widetilde L^2 (0,T ;  B^{\frac{1}{2}}_{\infty,\infty})}
\\
&
 \quad 
    + \| u_{m+1} \|_{\widetilde L^2 (0,T ;  B^{\frac{1}{2}}_{\infty,\infty})}
  \big) 
  \| u_n - u_m 
  \|_{\widetilde L^2 (0,T ;  B^{-\delta}_{\infty,q})}.
\end{split}
\end{gather}
\end{lem}

\noindent 
{\bf Remark}. 
The above lemma is concerned with the inhomogeneous Besov spaces, so that 
the constant $C$ depends on the time $T > 0$. 
On the other hand, we can also consider the homogeneous Besov spaces. 
In that case, the constant $C$ is independent of $T$. 

\vskip3mm

\begin{pf}
First we show the inequality \eqref{0402-3}. 
For $j = 5,6,7, \cdots$, it follows from the recurrence relation 
\eqref{0402-2} that 
\begin{equation}\label{0606-4}
\begin{split}
& 
\partial _t \phi_j * u_{n+1} + \Lambda \phi_j * u_{n+1} 
+   \big( S_{j-3} u_n\big) 
    \partial _x \phi_j * u_{n+1} 
\\
\displaystyle  
= 
& 
    \big( S_{j-3} u_n\big) 
    \partial _x \phi_j * u_{n+1} 
- \phi_j * \sum _{ l \geq 0 } 
    \big( S_{j-3} u_n\big) 
    \partial _x \phi_l * u_{n+1} 
\\
& -  \phi_j * \sum _{ l \geq 0} 
    \Big( \sum_{k \geq l+3} \phi_k * u _{n} 
    \Big) 
    \partial _x \phi_l * u_{n} 
  - \frac{1}{2} \phi_j * \partial _x \sum _{ |k-l| \leq 2} 
     (\phi_k * u_n) (\phi_l * u_n) .  
\end{split}
\end{equation}
By Proposition~\ref{prop:1}, we get that 
\begin{equation}\notag 
\begin{split}
& 
\partial _t \| \phi_j * u_{n+1} \|_{L^\infty} 
  + c 2^ j \| \phi_j * u_{n+1} \|_{L^\infty}
\\
\displaystyle  
\leq  
& 
    \Big\| \big( S_{j-3} u_n\big) 
    \partial _x \phi_j * u_{n+1} 
- \phi_j * \sum _{ l \geq 0} 
    \big( S_{j-3} u_n\big) 
    \partial _x \phi_l * u_{n+1}    \Big\|_{L^\infty}
\\
& + \Big\| \phi_j * \sum _{ l \in \mathbb Z} 
    \Big( \sum_{k \geq l+3} \phi_k * u _{n} 
    \Big) 
    \partial _x \phi_l * u_{n} 
    \Big\|_{L^\infty}
   + \frac{1}{2} \Big\| \phi_j * \partial _x \sum _{ |k-l| \leq 2} 
     (\phi_k * u_n) (\phi_l * u_n) 
\Big\|_{L^\infty}. 
\end{split}
\end{equation}
Noting that the left member is 
$e^{-t c2^{j}}\partial_t (e^{t c2^{j}} \|\phi_j * u_{n+1} \|_{L^\infty})$, 
multiplying by $e^{tc 2^j}$  
and integrating in both sides, we have that 
\begin{equation}\label{0606-2}
\| \phi_j * u_{n+1} (t) \|_{L^\infty} 
\leq e^{-tc 2^j} \| \phi_j * u_0\|_{L^\infty} 
+ {\rm I} (t) + {\rm II} (t) + {\rm III}(t) , 
\end{equation}
where 
\begin{equation}\notag 
\begin{split}
{\rm I} (t)
:= 
& \int_0^t e^{-(t-\tau)c 2^j}
 \Big\| \big( S_{j-3} u_n\big) 
    \partial _x \phi_j * u_{n+1} 
- \phi_j * \sum _{ l \geq 0} 
    \big( S_{j-3} u_n\big) 
    \partial _x \phi_l * u_{n+1}   
   \Big\|_{L^\infty}    d\tau ,
\\
{\rm II} (t)
:= 
& \int_0^t e^{-(t-\tau)c 2^j}
\Big\| \phi_j * \sum _{ l \in \mathbb Z} 
    \Big( \sum_{k \geq l+3} \phi_k * u _{n} 
    \Big) 
    \partial _x \phi_l * u_{n} 
    \Big\|_{L^\infty}  d\tau ,
\\
{\rm III} (t)
:= 
& \int_0^t e^{-(t-\tau)c 2^j}
\frac{1}{2} \Big\| \phi_j * \partial _x \sum _{ |k-l| \leq 2} 
     (\phi_k * u_n) (\phi_l * u_n) 
\Big\|_{L^\infty}  d\tau .
\end{split}
\end{equation}
We estimate the above three by the use of the Fourier mutiplier theorem and 
the H\"older inequality. 
The first term {\rm I} is estimated with a kind of commutator estimates 
(see e.g. (3.4) and (3.5) in \cite{Iw-2015}) as 
\begin{equation}\notag 
\begin{split}
{\rm I} (t) 
\leq 
& \int_0^t  \| \partial _x S_{j-3} u_n  \|_{L^\infty} 
   \sum _{\mu = -3}^3 \| \phi_{j+\mu} * u_{n+1} \|_{L^\infty}
   d\tau 
\\
\leq 
&C \int_0^t \sum _{ k \leq j-3} 2^k \| \phi_k * u_n \|_{L^\infty}
   \sum _{\mu = -3}^3 \| \phi_{j+\mu} * u_{n+1} \|_{L^\infty}
   d\tau 
\\
\leq 
& C 
 \sum _{ k \leq j-3} 2^k 
    \| \phi_k * u_n \|_{L^2 (0,T ; L^\infty)}
    \sum _{ \mu = -3}^3 
    \| \phi_{j+\mu} * u_{n+1} \|_{L^2 (0,T ; L^\infty)}
\\
\leq 
& C 
 \| u_n \|_{\widetilde L^2 (0,T ;  B^\frac{1}{2}_{\infty,\infty})}
 \sum _{ k \leq j-3} 2^{\frac{1}{2}k}
    \sum _{ \mu = -3}^3 
    \| \phi_{j+\mu} * u_{n+1} \|_{L^2 (0,T ; L^\infty)}. 
\end{split}
\end{equation}
We multiply by $2^{sj}$ and take the sequence norm $\ell ^q (\mathbb Z)$ to get 
\begin{equation}\notag 
\Big\{ \sum _{ j \geq 5} \big( 2^{sj}{\rm I } (t) \big)^q \Big\} ^{\frac{1}{q}} 
\leq C 
 \| u_n \|_{\widetilde L^2 (0,T ;  B^\frac{1}{2}_{\infty,\infty})}
 \| u_{n+1} \|_{\widetilde L^2 (0,T ;  B^{\frac{1}{2}+s}_{\infty,q})}.
\end{equation}
Since 
$
\phi_j * \sum _{ l \geq 0} 
    \big( \sum_{k \geq l+3} \phi_k * u _{n} 
    \big) 
    \partial _x \phi_l * u_{n} 
= 
\phi_j * 
\big( \sum _{ k = j-3 }^{j+3} 
    ( \phi_{k} * u _{n} )
    \sum _{ l \leq k -3} \partial _x \phi_l * u_{n} 
\big), 
$
the second term {\rm II} is also handled in the similar way 
to the first one ${\rm I}(t)$:  
\begin{equation}\notag 
\begin{split}
\Big( \sum _{ j \geq 5} \big( 2^{sj} {\rm II} (t) \big) ^q \Big\} ^{\frac{1}{q}}
\leq 
& C 
 \| u_n \|_{\widetilde L^2 (0,T ;  B^\frac{1}{2}_{\infty,\infty})}
 \| u_n \|_{\widetilde L^2 (0,T ;  B^{\frac{1}{2}+s}_{\infty,q})} .
\end{split}
\end{equation}
As to the third term {\rm III}, 
we also apply the Fourier multiplier theorem and the H\"older inequality 
to get that 
\begin{equation}\notag
\begin{split}
{\rm III} (t) 
\leq 
& C 
 \sum _{k \geq j-5} \sum _{l = k-2} ^{k+2}
 \int_0^t 2^j \| \phi_k * u_n \|_{L^\infty} \| \phi_l* u_n \|_{L^\infty}
 d\tau
\\
\leq 
& C 
 \sum _{k \geq j-5} 2^{-(k-j)} 
  2^k\sum _{l = k-2} ^{k+2}
 \| \phi_k * u_n \|_{L^2(0,T; L^\infty)} 
 \| \phi_l* u_n \|_{L^2(0,T ; L^\infty)}
\\
\leq 
& C 
 \sum _{m \geq -5} 2^{-m} 
 \cdot 
  2^{j+m}\sum _{\mu = -2} ^{2}
 \| \phi_{j+m} * u_n \|_{L^2(0,T; L^\infty)} 
 \| \phi_{j+m+\mu}* u_n \|_{L^2(0,T ; L^\infty)}. 
\end{split}
\end{equation}
By multiplying by $2^{sj}$ and taking the sequence norm of 
$\ell ^q (\mathbb Z)$, we obtain that 
\begin{equation}\notag 
\begin{split}
& 
\Big\{ \sum _{j \geq 5} 
  \big( 2^{sj}{\rm III} (t) \big) ^q \Big\} ^{\frac{1}{q}} 
\\
\leq 
& 
C \sum _{m \geq -5} 2^{-(1+s)m} 
\Big\{ \sum _{ j \geq 5} \Big(
  2^{ (1+s)(j+m)}\sum _{\mu = -2} ^{2}
\\
& \hspace*{35mm} 
\| \phi_{j+m} * u_n \|_{L^2(0,T; L^\infty)} 
 \| \phi_{j+m+\mu}* u_n \|_{L^2(0,T ; L^\infty)}
\Big) ^q \Big\}^{\frac{1}{q}}
\\
\leq 
& 
C  \| u_n \|_{\widetilde L^2 (0,T ;  B^\frac{1}{2}_{\infty,\infty})}
   \| u_n \|_{\widetilde L^2 (0,T ;  B^{\frac{1}{2}+s}_{\infty,q})} .
\end{split}
\end{equation}
By the above estimates for {\rm I, II, III} and the inequality \eqref{0606-2}, 
we obtain the inequality \eqref{0402-3} for the frequency away from the origin. 
As to the frequecy around the origin, 
we apply \eqref{0606-1}, \eqref{0606-1_inho} to \eqref{0606-4} 
and integrate to get that for $j = 0,1,2,3,4$ 
\begin{equation}\label{0608-2}
\begin{split}
\| \phi_j * u_{n+1} \| _{L^\infty} 
\leq 
& \| \phi_j * u_0 \|_{L^\infty} 
+ 
\int_0^t 
\Big( 
  \Big\| \phi_j * \sum _{ l \geq 0 } 
    \big( S_{j-3} u_n\big) 
    \partial _x \phi_l * u_{n+1} 
  \Big\|_{L^\infty}
\\
& \qquad 
+ \Big\|  \phi_j * \sum _{ l \geq 0} 
    \Big( \sum_{k \geq l+3} \phi_k * u _{n} 
    \Big) 
    \partial _x \phi_l * u_{n} 
    \Big\|_{L^\infty}
\\
&  \qquad  + \frac{1}{2} 
    \Big\| 
      \phi_j * \partial _x \sum _{ |k-l| \leq 2} 
     (\phi_k * u_n) (\phi_l * u_n) 
    \Big\|_{L^\infty}
    \Big) d\tau 
\\
\leq 
& \| \phi_j * u_0 \|_{L^\infty} 
+ 
C \int_0^t 
\Big( \sum _{ l = 0 }^5 \sum _{ k = 0}^2
    \| \phi_k * u_n \|_{L^\infty}
    \| \phi_l * u_{n+1} \|_{L^\infty}
\\
& \qquad 
+ 
 \sum_{k=0}^5 
    \sum_{l=0}^2 \| \phi_k * u _{n} \|_{L^\infty}
    \| \phi_l * u_{n} \|_{L^\infty}
\\
& \qquad 
   +
    \sum _{ |k-l| \leq 2} 2^{-(1+s)k}
      \cdot 
      2^{\frac{1}{2}k}  \| \phi_k * u_n\|_{L^\infty} 
      2^{(\frac{1}{2}+s)l} \| \phi_l * u_n) \|_{L^\infty}
    \Big) d\tau .
\end{split}
\end{equation}
Hence, we obtain that 
\begin{equation}\notag
\begin{split}
\Big\{ \sum _{ j = 0}^4 
  \Big( 2^{sj} \| \phi_j * u_{n+1} \|_{L^\infty (0,T ; L^\infty)}
  \Big) ^q 
 \Big\} ^{\frac{1}{q}}
 \leq 
& \| u_0 \|_{ B^{s}_{\infty,q}} 
+  C  \| u_n \|_{\widetilde L^2 (0,T ;  B^\frac{1}{2}_{\infty,\infty})}
      \| u_n \|_{\widetilde L^2 (0,T ;  B^{\frac{1}{2}+s}_{\infty,q})}
\\
&
+  C  \| u_n \|_{\widetilde L^2 (0,T ;  B^\frac{1}{2}_{\infty,\infty})}
      \| u_n \|_{\widetilde L^2 (0,T ;  B^{\frac{1}{2}+s}_{\infty,q})}, 
\end{split}
\end{equation}
which prove the estimate of \eqref{0402-3} for low frequency part. 
We completes the proof of \eqref{0402-3}.

We next prove the estimate \eqref{0402-4}. 
By taking $L^r (0,T)$ norm for the estimate \eqref{0606-2}, 
applying the inequality \eqref{0512-13} to the first term in the right member 
of \eqref{0606-2}
and order exchanging of integration for integrals of {\rm I, II} and {\rm III}, 
we have that for $j = 5,6,\cdots$
\begin{equation}\label{0907-1}
\| \phi_j * u_{n+1} (t) \|_{L^r (0,T ; L^\infty)} 
\leq C \Big( \| \phi_j * e^{-tc \Lambda} u_0 \|_{\widetilde L^r (0,T ;  L^\infty )} 
+  2^{-\frac{1}{r}j} 
 \big( \, \widetilde {\rm I} (T) + \widetilde {\rm II} (T) 
       + \widetilde {\rm III}(T) \big) 
\Big), 
\end{equation}
where $\widetilde {\rm I}, \widetilde {\rm II}, \widetilde {\rm III}$ 
are similar to ${\rm I}, {\rm II}, {\rm III}$ such that 
$e^{-(t-\tau)c 2^j}$ is replaced with $1$. 
By multiplying this inequality by $2^{(\frac{1}{r}+s)j}$, and 
taking the sequence norm of $\ell ^q (\mathbb Z)$, 
the same argument as before enables us to get the estimate of high 
frequency part of \eqref{0402-4}. 
As to the low frequency, we can also apply the argument \eqref{0608-2} 
to obtain the required estimate.

Let us prove the last estimate \eqref{0402-d}. 
By the recurrence relation \eqref{0402-2}, we write 
\begin{equation}\notag 
\begin{split}
\displaystyle 
& \partial _t ( u_{n+1} - u_{m+1} ) + \Lambda ( u_{n+1} - u_{m+1} )
+ \sum _{ l \in \mathbb Z} 
    \big( S_{ l-3} u _{n} \big) 
    \partial _x \phi_l * ( u_{n+1} - u_{m+1} )
\\
\displaystyle 
\quad 
=
& - \sum _{ l \in \mathbb Z} 
    \big( S_{ l-3} ( u _{n} - u_m) \big) 
    \partial _x \phi_l * u_{m+1} 
\\
& 
     -  \sum _{ l \in \mathbb Z} 
    \Big( \sum_{k \geq l+3} \phi_k * ( u _{n} - u_m )
    \Big) 
    \partial _x \phi_l * u_{n} 
     -  \sum _{ l \in \mathbb Z} 
    \Big( \sum_{k \geq l+3} \phi_k * u_m
    \Big) 
    \partial _x \phi_l * ( u_{n} - u_m )
\\
&  - \frac{1}{2} \partial _x \sum _{ |l-k| \leq 2} 
     \big(\phi_k * (u_n - u_m) \big) (\phi_j * u_n) 
   - \frac{1}{2} \partial _x \sum _{ |l-k| \leq 2} 
     (\phi_k * u_m ) \big( \phi_j * (u_n - u_m) \big) . 
\end{split}
\end{equation}
The similar arguments to the proof of \eqref{0402-3}, \eqref{0402-4} 
and Step 3 in the proof of Theorem~1.3 in \cite{WZ-2011} 
are applicable, each terms can be handled analogously to the previous 
estimates, 
and we obtain the estimate \eqref{0402-d}. 
\end{pf}


\vskip3mm

\noindent 
{\bf Proof of unique solvability in Theorem~\ref{thm:Burgers}.} 
First we derive a uniform boundedness 
of $\{ u_{n} \}_{n=1}^\infty$ to construct a solution. 
We consider the estimates that 
\begin{equation}\label{0608-4}
\| u_{n} \|_{\widetilde{L} ^\infty(0,T; B^{0}_{\infty,\infty})} 
\leq C_0 \| u_0 \|_{B^0_{\infty,\infty}}, 
\quad 
\| u_{n} \|_{\widetilde{L} ^2 (0,T; B^{\frac{1}{2}}_{\infty,\infty})} 
\leq 2\varepsilon,
\end{equation}
where the constant $C_0 \geq 2$ is larger than absolute constants 
appearing in the propositions and lemmas 
and $\varepsilon $ will be fixed as a small constant.  

When $n = 1$, \eqref{0608-4} is possible to be obtained, 
since the first one is just the boundedness of $e^{-t\Lambda}$ and 
the second one for small $T$ is assured by \eqref{0512-12} 
and the interpolation inequality, 
$$
\| f \|_{\widetilde L^2 (0,T ; B^{\frac{1}{2}}_{\infty,\infty} )}
\leq 
\| f \|_{\widetilde L^\infty (0,T ; B^{0}_{\infty,\infty} )} ^{\frac{1}{2}}
\| f \|_{\widetilde L^1 (0,T ; B^{1}_{\infty,\infty} )}^{\frac{1}{2}}, 
$$
which is on the controllability 
of the norm of $\widetilde L^2 (0,T ; B^{\frac{1}{2}}_{\infty,\infty} (\mathbb R))$ 
by $\widetilde L^\infty (0,T ; B^0_{\infty , \infty} (\mathbb R)) 
\cap \widetilde L^1 (0,T ; B^1_{\infty,\infty} (\mathbb R) )$. 
We also take $T$ smaller such that the inequality 
$
\| e^{-tc \Lambda} u_0 \|_{\widetilde L^2 (0,T ;  B^{\frac{1}{2}}_{\infty,\infty})}
\leq \varepsilon  
$ holds, 
where $c>0$ is a small constant appearing in \eqref{0402-4}.

Let us consider the estimates for $u_{n+1}$ under the assumption \eqref{0608-4} 
for $u_n$. 
It follows from the inequality \eqref{0402-3} for $s = 0$ 
and the assumption for $u_n$ that 
\begin{equation}\notag 
\begin{split}
\| u_{n+1} \|_{\widetilde{L} ^\infty(0,T; B^{0}_{\infty,\infty})} 
\leq  
&
\| u_0 \|_{ B^{0}_{\infty,\infty}} 
+  C_0  \frac{1}{16C_0} 
      \| u_{n+1} \|_{\widetilde L^2 (0,T ;  B^{\frac{1}{2}}_{\infty,\infty})}
+  C_0 \frac{1}{16C_0}
      \| u_n \|_{\widetilde L^2 (0,T ;  B^{\frac{1}{2}}_{\infty,\infty})}, 
\\
\| u_{n+1} \|_{\widetilde{L} ^\infty(0,T; B^{0}_{\infty,\infty})} 
\leq 
& 
2 \| u_0 \|_{B^0_{\infty,\infty}} .
\end{split}
\end{equation}
We also apply \eqref{0402-4} for $s = 0$
to get that 
\begin{equation}\notag 
\begin{split}
\| u_{n+1} \|_{\widetilde{L} ^2 (0,T; B^{\frac{1}{2}}_{\infty,\infty})} 
\leq  
&
\varepsilon  +  C_0  \cdot \frac{1}{16C_0}
      \| u_{n+1} \|_{\widetilde L^2 (0,T ;  B^{\frac{1}{2}}_{\infty,\infty})}
 +  C_0 \cdot \frac{1}{16 C_0}
      \| u_n \|_{\widetilde L^2 (0,T ;  B^{\frac{1}{2}-s}_{\infty,\infty})}, 
\\
\| u_{n+1} \|_{\widetilde{L} ^2 (0,T; B^{\frac{1}{2}}_{\infty,\infty})} 
\leq 
& 
2 \varepsilon .
\end{split}
\end{equation}
Hence, the uniform estimates \eqref{0608-4} is proved.

We next consider the convergence of $u_n$ in 
$\widetilde L^2 (0,T ;B^{-\delta}_{\infty,\infty} (\mathbb R))$. 
For fixed $0 < \delta < 1/2$,  
we have from the inequality \eqref{0402-d} that 
\begin{gather}\notag 
\begin{split}
\| u_{n+1} - u_{m+1} \|
_{\widetilde L^2 (0,T ;  B^{-\delta}_{\infty,\infty})}
\leq 
& 
C_0 \| e^{-tc \Lambda} (S_{n+1}u_0 - S_{m+1} u_0) 
  \|_{\widetilde L^2 (0,T ;  B^{-\delta}_{\infty,\infty})} 
\\
&
+ C_0 
 \cdot \frac{1}{16C_0} 
  \| u_{n+1} - u_{m+1} 
  \|_{\widetilde L^2 (0,T ;  B^{-\varepsilon }_{\infty,\infty} )}
\\
& 
+ C_0 \cdot \frac{3}{16C_0} 
  \| u_n - u_m 
  \|_{\widetilde L^2 (0,T ;  B^{-\delta}_{\infty,\infty})}, 
\end{split}
\\ \notag 
\begin{split}
\frac{15}{16} \| u_{n+1} - u_{m+1} \|
_{\widetilde L^2 (0,T ;  B^{-\delta}_{\infty,\infty})}
\leq 
& 
C \| (S_{n+1} - S_{m+1}) u_0 \|_{B^{-\frac{1}{2}-\delta}_{\infty,\infty}} 
\\
& 
+ \frac{3}{16} 
  \| u_n - u_m 
  \|_{\widetilde L^2 (0,T ;  B^{-\delta}_{\infty,\infty})},
\end{split}
\\ \notag 
\begin{split}
\| u_{n+1} - u_{m+1} \|
_{\widetilde L^2 (0,T ;  B^{-\delta}_{\infty,\infty})}
\leq 
& 
C 2^{- \min\{ n , m \}} 
  \| u_0 \|_{B^0_{\infty,\infty}}
+ \frac{1}{5} 
  \| u_n - u_m 
  \|_{\widetilde L^2 (0,T ;  B^{-\delta}_{\infty,\infty})}.
\end{split}
\end{gather}
The above inequality yeilds that 
$$
\| u_{n+1} - u_{n} \|
_{\widetilde L^2 (0,T ;  B^{-\delta}_{\infty,\infty})} 
\leq C 2^{-n} \| u_0 \|_{B^0_{\infty,\infty}},
$$
which prove the existence of the following limit 
$$
u:=\lim _{n \to \infty} u_n 
= \lim_{n\to \infty} 
    u_1 +\sum _{k=1} ^{n-1} (u_{k+1} - u_{k})
    \quad 
\text{in } \widetilde L^2 (0,T ;  B^{-\delta}_{\infty,\infty} (\mathbb R).
$$
It can be checked that $u$ also satisfies the same inequality 
as \eqref{0608-4}, 
and the standard limit argument ensures that $u$ 
satisfies \eqref{eq:B} in the sense of distribution. 
The uniqueness is proved by applying the inequality like \eqref{0402-d} 
for the same initial data. 
\hfill $\Box$

\vskip3mm

We next study the analyticity of solutions. 
Let us prepare a lemma and a proposition.

\begin{prop}\label{prop:2}
Let $\alpha, \beta \in \mathbb N \cup \{ 0 \}$ and 
$1 \leq q \leq \infty$. 
Then there exists a positive constant $C_0$ 
independent of $\alpha , \beta$ such that 
for any $u_0 \in  B^0_{\infty,q} (\mathbb R)$  
\begin{equation}\label{0606-3}
\| t^{\alpha + \beta} \partial_t ^\alpha \partial_x ^\beta e^{-t\Lambda} u_0 
\|_{B^0_{\infty,q} }
+
\| t^{\alpha + \beta} \partial_t ^\alpha \partial_x ^\beta e^{-t\Lambda} u_0 
\|_{\widetilde L^1 (0,T ; B^1_{\infty,q}) }
  \leq 
   C_0^{\alpha + \beta} (\alpha + \beta) ! 
     \, \| u_0 \|_{ B^0_{\infty,q}}. 
\end{equation}
\end{prop}

\begin{pf}
It follows from 
$\partial_t e^{-t\Lambda} = - \Lambda e^{-t\Lambda}$ that 
$$
\| \partial_t ^\alpha \partial_x ^\beta e^{-t\Lambda} u_0 
\|_{ B^0_{\infty,q} }
= 
  \| \Lambda ^{\alpha} \partial_x^\beta e^{-t\Lambda} u_0 
  \|_{ B^0_{\infty,q}} 
= 
  \| ( \Lambda  e^{-\frac{t}{\alpha + \beta}\Lambda} )^{\alpha} 
     ( \partial_x  e^{-\frac{t}{\alpha + \beta}\Lambda} )^{\beta} u_0 
  \|_{ B^0_{\infty,q}}.
$$
The smoothing effect of $e^{-t\Lambda}$ implies that 
$$
\| \partial_t ^\alpha \partial_x ^\beta e^{-t\Lambda} u_0 
\|_{ B^0_{\infty,q} }
\leq C ^{\alpha + \beta} 
 \Big( \frac{\alpha + \beta}{t} \Big) ^{\alpha + \beta} 
 \| u_0 \|_{ B^0_{\infty,q}}. 
$$
The above inequality and Stirling's approximation 
yield the desired estimate of the first term in the left member of 
\eqref{0606-3}. As for the second term, we also have from 
the above estimate and the maximal regularity estimate in the Chemin-Lerner spaces 
that 
\begin{equation}\notag
\begin{split}
\| t^{\alpha + \beta} \partial_t ^\alpha \partial_x ^\beta e^{-t\Lambda} u_0 
\|_{\widetilde L^1 (0,T ; B^1_{\infty,q})}
=
& 
\| ( \Lambda  e^{-\frac{t}{2( \alpha + \beta )}\Lambda} )^{\alpha} 
     ( \partial_x  e^{-\frac{t}{2(\alpha + \beta)}\Lambda} )^{\beta} 
     e^{\frac{t}{2}\Lambda}u_0 
\|_{\widetilde L^1 (0,T ; B^1_{\infty,q})}
\\
\leq 
& 
C ^{\alpha + \beta} ( \alpha + \beta) ^{\alpha + \beta} 
\|   e^{\frac{t}{2}\Lambda}u_0 
\|_{\widetilde L^1 (0,T ; B^1_{\infty,q})}
\\
\leq 
& 
C C ^{\alpha + \beta} ( \alpha + \beta) ^{\alpha + \beta} 
\| u_0 \|_{ B^1_{\infty,q} }.
\end{split}
\end{equation}
Therefore we obtain the inequality \eqref{0606-3}. 
\end{pf}

\vskip3mm 

Based on the proof of Lemma~\ref{lem:3}, 
we show the following nonlinear estimates to obtain the analyticity.

\begin{prop}\label{prop:anal_nonlinear}
Let 
$
\partial_{t,x} ^\alpha
= \partial_t ^{\alpha _0} \partial_{x}^{\alpha_1}
$, 
$\alpha ! := \alpha _0 ! \alpha _1 !$ 
and $|\alpha| = \alpha _0 + \alpha _1$ 
for $\alpha = (\alpha_0, \alpha_1)
\in (\mathbb N \cup \{ 0 \})^{2}$. 
Assume that $|\alpha | \geq 1$. 
Then there exists a positive constant $C$ independent of $\alpha$ such that 
for any $1 \leq q \leq \infty$ 
\begin{gather}\label{0606-5}
\begin{split}
& 
\| t^{|\alpha|} \partial_{t,x} ^\alpha u_{n+1} 
\|_{\widetilde L^\infty (0,T ;  B^0_{\infty,q}) 
  \cap \widetilde L^1 (0,T ; B^1_{\infty,q})} 
\\
\leq 
& C |\alpha| 
  \| t ^{|\alpha| -1} \partial_{t,x} ^\alpha u_{n+1} 
  \|_{\widetilde L^1 (0,T ; B^0_{\infty,q})} 
\\
&+ C \| u_n \|_{\widetilde L^1 (0,T ;  B^1_{\infty,1})}
    \| t^{|\alpha|}\partial_{t,x} ^\alpha u_{n+1} 
    \|_{\widetilde L^\infty (0,T ;  B^0_{\infty,q})}
\\
&+ C \sum _{\beta \not = 0 , \, \beta + \gamma = \alpha} 
    \frac{\alpha !}{\beta ! \gamma !} 
    \| t^{|\beta|} \partial _{t,x} ^\beta u_n 
    \|_{\widetilde L^\infty (0,T ;  B^0_{\infty,1})}
    \| t^{|\gamma|}\partial_{t,x} ^\gamma u_{n+1} 
    \|_{\widetilde L^1 (0,T ;  B^1_{\infty,q})}
\\
& + C\sum _{\beta + \gamma = \alpha} 
    \frac{\alpha !}{\beta ! \gamma !} 
    \| t^{|\beta|} \partial _{t,x} ^\beta u_n 
    \|_{\widetilde L^2 (0,T ;  B^\frac{1}{2}_{\infty,\infty})} 
    \| t^{|\gamma|} \partial _{t,x} ^\gamma u_n 
    \|_{\widetilde L^2 (0,T ;  B^\frac{1}{2}_{\infty,q})} ,
\end{split}
\end{gather}
\end{prop}

\begin{pf}
By applying $\partial_{t,x} ^\alpha $ to 
the recurrence relation \eqref{0402-2} and the similar argument 
as in \eqref{0606-4}, we write 
\begin{equation}\notag 
\begin{split}
& 
\partial _t \phi_j * ( \partial_{t,x} ^\alpha u_{n+1} )
+ \Lambda \phi_j * (\partial_{t,x} ^\alpha u_{n+1} )
+   \big( S_{j-3} u_n \big) 
    \partial _x \phi_j * (\partial_{t,x} ^\alpha u_{n+1} )
\\
\displaystyle  
= 
& 
    \big( S_{j-3} u_n\big) 
    \partial _x \phi_j * (\partial_{t,x} ^\alpha u_{n+1} )
- \phi_j * \sum _{ l \in \mathbb Z} 
    \big( S_{j-3} u_n\big) 
    \partial _x \phi_l * (\partial_{t,x} ^\alpha u_{n+1} )
\\
& 
+ \phi_j * \sum _{ l \in \mathbb Z} 
    \big( S_{j-3} u_n\big) 
    \partial _x \phi_l * (\partial_{t,x} ^\alpha u_{n+1} ) 
- 
\partial_{t,x} ^\alpha \phi_j * \sum _{ l \in \mathbb Z} 
    \big( S_{j-3} u_n\big) 
    \partial _x \phi_l * u_{n+1} 
\\
& - \partial_{t,x} ^\alpha  \phi_j * \sum _{ l \in \mathbb Z} 
    \Big( \sum_{k \geq l+3} \phi_k * u _{n} 
    \Big) 
    \partial _x \phi_l * u_{n} 
  - \partial_{t,x} ^\alpha 
    \frac{1}{2} \phi_j * \partial _x \sum _{ |k-l| \leq 2} 
     (\phi_k * u_n) (\phi_l * u_n) .  
\end{split}
\end{equation}
It follows from Proposition~\ref{prop:1} and the Leipniz rule that 
for $j = 5,6, \cdots$ 
\begin{equation}\notag 
\begin{split}
& 
\partial_t \| \phi_j * (\partial_{t,x} ^\alpha u_{n+1})\|_{L^\infty} 
+ c 2^ j \| \phi_j * (\partial_{t,x} ^\alpha u_{n+1})\|_{L^\infty} 
\\
\leq 
& \Big\| \big( S_{j-3} u_n\big) 
    \partial _x \phi_j * (\partial_{t,x} ^\alpha u_{n+1} )
- \phi_j * \sum _{ l \in \mathbb Z} 
    \big( S_{j-3} u_n\big) 
    \partial _x \phi_l * (\partial_{t,x} ^\alpha u_{n+1} )
\Big\| _{L^\infty}
\\
& 
+   \sum _{\beta \not = 0, \, \beta + \gamma = \alpha} 
    \frac{\alpha !}{\beta ! \gamma !} 
  \Big\|
    \phi_j * \sum _{ l \in \mathbb Z} 
    \big( S_{j-3} \partial _{t,x} ^\beta u_n\big) 
    \partial _x \phi_l * ( \partial _{t,x} ^\gamma u_{n+1} )
  \Big\|_{L^\infty}
\\
& +\sum _{\beta + \gamma = \alpha} 
    \frac{\alpha !}{\beta ! \gamma !} 
  \Big\| 
    \phi_j * \sum _{ l \in \mathbb Z} 
    \Big( \sum_{k \geq l+3} \phi_k * \partial_{t,x}^\beta u _{n} 
    \Big) 
    \partial _x \phi_l * ( \partial _{t,x} ^\gamma u _{n} \big) 
  \Big\|_{L^\infty}
\\
& + \frac{1}{2}\sum _{\beta + \gamma = \alpha} 
    \frac{\alpha !}{\beta ! \gamma !} 
  \Big\| 
     \phi_j * \partial _x \sum _{ |k-l| \leq 2} 
     (\phi_k * \partial_{t,x}^\beta u_n) (\phi_l * \partial _{t,x} ^\gamma u_n) 
  \Big\|_{L^\infty} . 
\end{split}
\end{equation}
Multiplying the above inequality by $t^{|\alpha|}$ and noting that 
$t^{|\alpha|} \partial _t f = \partial_t (t^{|\alpha|} f) - |\alpha| t^{|\alpha|  -1} f$ 
with $f = \| \phi_j * (\partial _{t,x} ^\alpha u_{n+1}) \|_{L^\infty}$ 
and $t^{|\alpha|} = t ^{|\beta|} t^{|\gamma|}$, 
we have that 
\begin{equation}\notag 
\begin{split}
& 
\partial_t \| \phi_j * (t^{|\alpha|}\partial_{t,x} ^\alpha u_{n+1})\|_{L^\infty} 
+ c 2^ j \| \phi_j * (t^{|\alpha|}\partial_{t,x} ^\alpha u_{n+1})\|_{L^\infty} 
\\
\leq 
& |\alpha| 
  \| \phi_j * (t ^{|\alpha| -1} \partial_{t,x} ^\alpha u_{n+1})\|_{L^\infty} 
\\
& 
+ \Big\| \big( S_{j-3} u_n\big) 
    \partial _x \phi_j * ( t ^{|\alpha|}\partial_{t,x} ^\alpha u_{n+1} )
- \phi_j * \sum _{ l \in \mathbb Z} 
    \big( S_{j-3} u_n\big) 
    \partial _x \phi_l * (t ^{|\alpha|} \partial_{t,x} ^\alpha u_{n+1} )
\Big\| _{L^\infty}
\\
& 
+   \sum _{\beta \not = 0, \, \beta + \gamma = \alpha } 
    \frac{\alpha !}{\beta ! \gamma !} 
  \Big\|
    \phi_j * \sum _{ l \in \mathbb Z} 
    \big( S_{j-3} t ^{|\beta|}\partial _{t,x} ^\beta u_n\big) 
    \partial _x \phi_l * ( t ^{|\gamma|}\partial _{t,x} ^\gamma u_{n+1} )
  \Big\|_{L^\infty}
\\
& +\sum _{\beta + \gamma = \alpha} 
    \frac{\alpha !}{\beta ! \gamma !} 
  \Big\| 
    \phi_j * \sum _{ l \in \mathbb Z} 
    \Big( \sum_{k \geq l+3} \phi_k * t ^{|\beta|} \partial_{t,x}^\beta u _{n} 
    \Big) 
    \partial _x \phi_l * ( t ^{|\gamma|} \partial _{t,x} ^\gamma u _{n} \big) 
  \Big\|_{L^\infty}
\\
& + \frac{1}{2}\sum _{\beta + \gamma = \alpha} 
    \frac{\alpha !}{\beta ! \gamma !} 
  \Big\| 
     \phi_j * \partial _x \sum _{ |k-l| \leq 2} 
     (\phi_k * t ^{|\beta|} \partial_{t,x}^\beta u_n) 
     (\phi_l * t ^{|\gamma|} \partial _{t,x} ^\gamma u_n) 
  \Big\|_{L^\infty} . 
\end{split}
\end{equation}
By $t^{|\alpha|} \partial_{t,x} ^\alpha  u_{n+1} = 0$ for $t = 0$ 
and 
the analogous argument to the estimates of {\rm I, II, III} 
and ${\rm \widetilde I, \widetilde{II}, \widetilde {III} } $
appearing in \eqref{0606-2} and \eqref{0907-1}, 
we obtain the desired inequality \eqref{0606-5} 
for the frequency away from the origin. 
As for the frequency around the origin, we apply the similar inequality 
to \eqref{0608-2} to get the required inequality. 
The proof is finished. 
\end{pf}

\vskip3mm 

Let us turn to prove the analyticity. 

\vskip3mm 

\noindent 
{\bf Proof of analyticity in Theorem~\ref{thm:Burgers}.} 
The solution $u$ constructed in the previous proof of unique solvability 
satisfies 
$u (t) \in B^0_{\infty,1} (\mathbb R)$  for almost every $t \in (0,T]$ 
since 
$\widetilde L^1 (0,T ; B^1_{\infty,\infty} (\mathbb R))$ is 
embedded to $L^1 (0,T ; B^0_{\infty,1} (\mathbb R))$. 
Hence, the problem is reduced to the one for initial data 
$u_0 \in B^0_{\infty,1} (\mathbb R)$. 

For initial data $u_0 \in B^0_{\infty,1} (\mathbb R)$, 
we consider the sequence $\{ u_n \}_{n=1}^\infty $ in \eqref{0402-2} again 
to obtain the uniform boundedness in 
$\widetilde L^\infty (0,T ; B^0_{\infty,1} (\mathbb R)) 
\cap \widetilde L^1 (0,T ; B^1_{\infty,1} (\mathbb R))$ 
like \eqref{0608-4}: 
\begin{equation}\label{0608-6}
\| u_{n} \|_{\widetilde{L} ^\infty([0,T]; B^{0}_{\infty, 1})} 
\leq C_0 \| u_0 \|_{B^0_{\infty, 1}}, 
\quad 
\| u_{n} \|_{\widetilde{L} ^1 (0,T; B^{1}_{\infty,  1})} 
\leq 2\varepsilon,
\end{equation}
with small fixed constant $\varepsilon >0$. 
By Lemma~\ref{lem:3} and 
the similar argument to unique solvability, 
we obtain the uniform boundedness of $\{ u_n \}_{n=1}^\infty$ in 
the space 
$\widetilde L^\infty (0,T ; B^0_{\infty,q} (\mathbb R)) 
\cap \widetilde L^1 (0,T ; B^1_{\infty,q} (\mathbb R))$ 
for $q = 1,\infty$. In addition, we also derive 
another boundedness for analyticity, namely, we will show that 
there exist $C_0, C_1 > 0$ such that 
\begin{equation}\label{0608-7}
\| t^{|\alpha|} \partial_{t,x} ^\alpha u_{n} 
\|_{\widetilde L^\infty (0,T ;  B^0_{\infty,1})
    \cap \widetilde L^1 (0,T ; B^1_{\infty,1}) }
 \leq \frac{C_0 ^{|\alpha|-1} C_1 ^{|\alpha|} \alpha !}{(1+|\alpha|)^4} 
\end{equation}
for any $n \in \mathbb N$, 
$\alpha = (\alpha _0 , \alpha_1) \in (\mathbb N \cup \{ 0 \})^2$, 
where 
$\partial_{t,x} ^\alpha = \partial_t ^{\alpha_0} \partial_x ^{\alpha _1}$, 
$C_0$ and $C_1$ 
are fixed sufficiently large and 
depend on $\| u_0 \|_{B^0_{\infty,1}}$.

We prove \eqref{0608-7} by induction argument with respect to 
$n$ and $\alpha$.  
When $n = 1$, \eqref{0608-7} is proved by Proposition~\ref{prop:2}. 
When $|\alpha| = 0$, the inequality \eqref{0608-7} is true by \eqref{0608-6}. 
Here take $n' \in \mathbb N$ and 
$\alpha' \in (\mathbb N \cup \{ 0 \})^2 \setminus \{ (0,0) \}$ and 
let us assume that \eqref{0608-7} is true for $u_{n'}$ 
with all $\alpha \in (\mathbb N \cup \{ 0 \})^2$ 
and for $u_{n'+1}$ with all $\alpha$ such that $| \alpha| \leq |\alpha '| -1$. 
We need to consider the estimate \eqref{0608-7} for $u_{n' +1}$ with $\alpha = \alpha '$. 
It follows from the inequality \eqref{0606-5} and the assumption of induction 
and \eqref{0608-6} that 
\begin{equation}\notag 
\begin{split}
& 
\| t^{|\alpha'|} \partial_{t,x} ^{\alpha'} u_{n'+1} 
\|_{\widetilde L^\infty (0,T ;  B^0_{\infty,1}) 
  \cap \widetilde L^1 (0,T ; B^1_{\infty,1})} 
\\
\leq 
& C |\alpha'| 
  \| t ^{|\alpha'| -1} \partial_{t,x} ^{\alpha '} u_{n'+1} 
  \|_{\widetilde L^1 (0,T ; B^0_{\infty,1})} 
+ C \| u_{n'} \|_{\widetilde L^1 (0,T; B^1_{\infty,1})} 
     \| t^{|\alpha'|} \partial_{t,x}^{\alpha'} u_{n'+1} 
     \|_{\widetilde L^\infty (0,T ; B^0_{\infty,1})}
 \\
&+  C \sum _{\beta \not= 0, \beta + \gamma = \alpha'} 
    \frac{\alpha' !}{\beta ! \gamma !} 
    \frac{C_0 ^{|\beta|-1} C_1 ^{|\beta|} \beta !}{(1+|\beta|)^4} 
    \frac{C_0 ^{|\gamma|-1} C_1 ^{|\gamma|} \gamma !}{(1+|\gamma|)^4} 
\\
& + C \sum _{\beta + \gamma = \alpha'} 
    \frac{\alpha ' !}{\beta ! \gamma !} 
    \frac{C_0 ^{|\beta|-1} C_1 ^{|\beta|} \beta !}{(1+|\beta|)^4} 
    \frac{C_0 ^{|\gamma|-1} C_1 ^{|\gamma|} \gamma !}{(1+|\gamma|)^4} .
\\
\leq 
& C |\alpha'| 
  \| t ^{|\alpha'| -1} \partial_{t,x} ^{\alpha '} u_{n'+1} 
  \|_{\widetilde L^1 (0,T ; B^0_{\infty,1})} 
+ C \cdot 2 \varepsilon 
     \| t^{|\alpha '|} \partial_{t,x}^{\alpha '}u_{n'+1} 
     \|_{\widetilde L^\infty (0,T ; B^0_{\infty,1})}
\\
& + \frac{C C_0^{|\alpha'|-2} C_1 ^{|\alpha'|} \alpha ' !}{(1+|\alpha '|)^4},
\end{split}
\end{equation}
and 
\begin{equation}\notag 
\begin{split}
& 
\| t^{|\alpha'|} \partial_{t,x} ^{\alpha '} u_{n'+1} 
\|_{\widetilde L^\infty (0,T ;  B^0_{\infty,1}) 
  \cap \widetilde L^1 (0,T ; B^1_{\infty,1})} 
\\
\leq 
& C |\alpha '| 
  \| t ^{|\alpha '| -1} \partial_{t,x} ^\alpha u_{n'+1} 
  \|_{\widetilde L^1 (0,T ; B^0_{\infty,1})} 
+ \frac{C}{C_0} \cdot 
  \frac{ C_0^{|\alpha '|-1} C_1 ^{|\alpha '|} \alpha ' ! }{ (1+|\alpha '|)^4}, 
\end{split}
\end{equation}
since we can take $\varepsilon$ sufficiently small. Hence, we need to 
estimate 
$|\alpha '| 
  \| t ^{|\alpha '| -1} \partial_{t,x} ^{\alpha '} u_{n' +1} 
  \|_{\widetilde L^1 (0,T ; B^0_{\infty,q})} $. 
For $ \alpha '= (\alpha '_0 , \alpha ' _1 )$, 
we have 
\begin{equation}\notag 
\begin{split}
& 
|\alpha'| 
\| t^{|\alpha'| -1} \partial_{t,x}^{\alpha '} u_{n'+1}  
  \|_{L^1 (0,T ; B^0_{\infty,1})}
\\
= 
& 
\alpha '_0 
\| t^{\alpha '_0 -1 + \alpha '_1} \partial_t 
   (\partial_t ^{\alpha '_0-1} \partial_x ^{\alpha '_1}) u_{n' +1}  
  \|_{L^1 (0,T ; B^0_{\infty,1})}
+ 
\alpha '_1
\| t^{\alpha '_0 + \alpha '_1 -1} \partial_x 
   (\partial_t ^{\alpha '_0} \partial_x ^{\alpha '_1 -1}) u_{n'+1}  
  \|_{L^1 (0,T ; B^0_{\infty,1})} . 
\end{split}
\end{equation}
For the first term in the right member, 
by operating $t^{\alpha'_0 - 1 + \alpha'_1} 
\partial_t ^{\alpha' _0 - 1} \partial_x ^{\alpha' _1}$ to \eqref{0402-2},  
estimating directly with taking the norm of $L^1 (0,T ; B^0_{\infty,1} (\mathbb R))$, 
and applying the assumption of the induction for $\alpha $ with 
$| \alpha| \leq (\alpha'_0 -1) + \alpha'_1 = |\alpha' | -1$,  
we get that 
\begin{equation}\notag 
\begin{split}
& 
\alpha' _0 
\| t^{\alpha'_0 -1 + \alpha'_1} \partial_t 
   (\partial_t ^{\alpha'_0 -1} \partial_x ^{\alpha'_1} ) u_{n' +1}  
  \|_{L^1 (0,T ; B^0_{\infty,1})}
\\ 
\leq 
& 
C \alpha' _0  
\| t^{\alpha'_0 -1 + \alpha'_1} 
   \partial_t ^{\alpha'_0 -1} \partial_x ^{\alpha'_1} u_{n' +1}  
  \|_{L^1 (0,T ; B^1_{\infty,1})}
\\
& + C \sum _{\beta + \gamma \leq (\alpha'_0 -1 , \alpha'_1)} 
   \frac{(\alpha'_0 -1)!  \alpha'_1 !}{\beta ! \gamma !}
   \|t^{|\beta|} \partial_{t,x}^\beta u_{n'} 
   \|_{\widetilde L^\infty (0,T ; B^0_{\infty,1})}
   \| t^{|\gamma|} \partial_{t,x}^{\gamma} u_{n' +1}
   \|_{\widetilde L^1 (0,T ; B^1_{\infty,1})}
\\
& + C \sum _{\beta + \gamma \leq (\alpha'_0 -1 , \alpha'_1)} 
   \frac{(\alpha'_0 - 1)! \alpha'_1 !}{\beta ! \gamma !}
   \|t^{|\beta|} \partial_{t,x}^\beta u_{n'} 
   \|_{\widetilde L^\infty (0,T ; B^0_{\infty,1})}
   \| t^{|\gamma|} \partial_{t,x}^{\gamma} u_{n'}
   \|_{\widetilde L^1 (0,T ; B^1_{\infty,1})}
\\ 
\leq 
& 
C \alpha'_0 \cdot \frac{C_0^{|\alpha'| -2} C_1^{|\alpha'| -1} 
  (\alpha'_0 -1)! \alpha'_1  !}{(1+|\alpha'| -1)^4}
+ \frac{C}{C_0} \cdot 
  \frac{ C_0^{|\alpha'|-2} C_1 ^{|\alpha'|-1} \alpha' ! }{ (1+|\alpha'|)^4} 
\\ 
= 
& 
 \frac{2C}{C_0^2 C_1} \cdot 
  \frac{ C_0^{|\alpha'|-1} C_1 ^{|\alpha'|} \alpha' ! }{ (1+|\alpha'|)^4} .
\end{split}
\end{equation}
As for the second, since $\partial_x$ is a mapping from 
$B^1_{\infty,1} (\mathbb R)$ to $B^0_{\infty,1} (\mathbb R)$, 
the following holds: 
\begin{equation}\notag 
\begin{split}
\alpha'_1
\| t^{\alpha'_0 + \alpha'_1 -1} \partial_x 
   (\partial_t ^{\alpha'_0 } \partial_x ^{\alpha'_1 - 1}) u_{n'+1}  
  \|_{L^1 (0,T ; B^0_{\infty,1})}
\leq 
& C \alpha'_1 
\| t^{\alpha'_0 + \alpha'_1 -1} 
   \partial_t ^{\alpha'_0} \partial_x ^{\alpha'_1 - 1} u_{n' +1}  
  \|_{L^1 (0,T ; B^1_{\infty,1})}
\\
\leq 
& C \alpha'_1 \cdot 
 \frac{C_0^{|\alpha'|-2} C_1 ^{|\alpha'|-1} \alpha'_0 ! (\alpha'_1-1)!} 
      {(1+|\alpha'|-1)^4}
\\
= 
& \frac{C}{C_0 C_1} \cdot 
 \frac{C_0^{|\alpha'|-1} C_1 ^{|\alpha'|} \alpha' !} {(1+|\alpha'|)^4} . 
\end{split}
\end{equation}
Hence, the above four inequalities yield that 
\begin{equation}\notag 
\begin{split}
\| t^{|\alpha'|} \partial_{t,x} ^{\alpha'} u_{n'+1} 
\|_{\widetilde L^\infty (0,T ;  B^0_{\infty,1}) 
  \cap \widetilde L^1 (0,T ; B^1_{\infty,1})} 
\leq 
& 
\Big( \frac{C}{C_0} + \frac{C}{C_0^2 C_1} \Big) 
\frac{C_0^{|\alpha'|-1} C_1 ^{|\alpha'|} \alpha' !} {(1+|\alpha'|)^4} 
\\
\leq 
& 
\frac{C_0^{|\alpha'|-1} C_1 ^{|\alpha'|} \alpha' !} {(1+|\alpha'|)^4} ,
\end{split}
\end{equation}
where $C_0, C_1$ are taken sufficiently large. 
Therefore, we obtain the estimate \eqref{0608-7} for all 
$n \in \mathbb N$ and $\alpha \in ( \mathbb N \cup \{ 0 \} ) ^2$. 

The estimates \eqref{0608-6} implies the solvability of solution 
in the space 
$C([0,T], B^0_{\infty,1} (\mathbb R)) \cap 
\widetilde L^1 (0,T ; B^1_{\infty,1} (\mathbb R))$, 
and \eqref{0608-7} for all $n \in \mathbb N$ 
enables us to get the same boundedness of the solution,   
which verifies the analyticity of solution in space-time. 
The proof of analyticity is completed. 
\hfill $\Box$

\section{Large time behavior for Burgers equation}

We start by considering the initial data 
$u_0 \in L^1 (\mathbb R) \cap B^0_{\infty,1} (\mathbb R)$, 
since initial data in $L^1 (\mathbb R)$ with \eqref{B_initial}
can be approximated by functions 
in $L^1 (\mathbb R) \cap B^0_{\infty,1} (\mathbb R)$.  
Let $u_0 \in L^1 (\mathbb R) \cap B^0_{\infty,1} (\mathbb R)$ 
and we consider a solution $u$ of \eqref{eq:B} such that 
\begin{equation}\label{0723-5}
u \in C([0,\infty), B^0_{\infty,1} (\mathbb R)) 
\cap L^1_{loc} ([0, \infty), B^1_{\infty,1} (\mathbb R)),
\end{equation}
noting that the argument of solvability in section \ref{sec:LWP} is applicable 
to that in $B^0_{\infty,1} (\mathbb R)$ and 
the global regularity (see e.g. \cite{KNS-2008,MW-2009}) assures the 
global existence. 
We prepare the following lemma to guarantee the integrability of 
the solution for all the positive time.

\begin{lem}\label{lem:0723-1}
Let $u_0 \in L^1 (\mathbb R) \cap B^0_{\infty,1} (\mathbb R)$. 
Then there exists a unique global solution 
$u \in C([0,\infty), L^1 (\mathbb R) \cap B^0_{\infty,1} (\mathbb R)) 
\cap L^1_{loc} (0,\infty ; B^1_{\infty,1} (\mathbb R))$ 
of the integral equation 
\begin{equation}\label{0723-6}
u(t) = e^{-t\Lambda} u_0 
  - \frac{1}{2} \int_0^t e^{-(t-\tau)\Lambda} \partial_x u (\tau)^2 \, d\tau .
\end{equation}
\end{lem}

\begin{pf}
We start by the local solvability in $L^1 (\mathbb R) \cap B^0_{\infty,1} (\mathbb R)$. 
We consider an sequence approximating solutions in the framework of 
$B^1_{\infty,1} (\mathbb R)$ by the similar argument of 
the local solvability in Section~\ref{sec:LWP} 
together with the boundedness in $L^1 (\mathbb R)$. 
Let $\{  u_n \}_{n=1}^\infty$ be defined by \eqref{0402-2}. 
On the estimates in $L^1 (\mathbb R)$, we estimate 
the corresponding integral equation in \eqref{0402-2} that 
\begin{gather}\notag 
\| u_1 (t) \|_{L^1} = \| e^{-t\Lambda}u_0 \|_{L^1} \leq \| u_0 \|_{L^1}, 
\\ \notag 
\begin{split}
\| u_{n+1} (t) \|_{L^1} 
\leq 
& 
\| u_0 \|_{L^1} 
+ C \int_0^t \| u_n \|_{L^1} 
   (\| u_{n+1} \|_{\dot B^1_{\infty,1}} + \| u_n \|_{\dot B^1_{\infty,1}}) d\tau 
\\
\leq  
& 
\|u_0\|_{L^1} 
+ C (\| u_{n+1} \|_{L^1 (0,t ; \dot B^1_{\infty,1})} 
     + \| u_{n} \|_{L^1 (0,t ; \dot B^1_{\infty,1})}
    ) \| u \|_{L^\infty (0,t ; L^1)} .
\end{split}
\end{gather}
Following the proof of unique solvability in section~\ref{sec:LWP}, 
we have from Lemma~\ref{lem:3} for $q = 1$ 
that 
\begin{gather}\notag 
\| u_{n} \|_{\widetilde{L} ^\infty(0,T; B^{0}_{\infty,1})} 
\leq C \| u_0 \|_{B^0_{\infty,1}}, 
\quad 
\| u_{n} \|_{\widetilde{L} ^1 (0,T; B^{1}_{\infty,1})} 
\leq 2\varepsilon,
\\ \notag 
\begin{split}
\| u_{n+1} - u_{m+1} \|
_{\widetilde L^2 (0,T ;  B^{-\delta}_{\infty,1})}
\leq 
& 
C 2^{- \min\{ n , m \}} 
  \| u_0 \|_{B^0_{\infty,1}}
+ \frac{1}{5} 
  \| u_n - u_m 
  \|_{\widetilde L^2 (0,T ;  B^{-\delta}_{\infty,1})}, 
\end{split}
\end{gather}
where $\varepsilon >0$ is a fixed sufficiently small constant, 
$T$ is small and is depending on $u_0$. 
Hence we deduce from the above four estimates that 
for some fixed small constant $\varepsilon > 0$ 
and small $T$, $u_n$ satisfies that 
\begin{gather}\notag 
\| u_{n} \|_{L ^\infty(0,T; L^1)} \leq C \| u_0 \|_{L^1},
 \quad 
\| u_{n} \|_{\widetilde{L} ^\infty(0,T; B^{0}_{\infty,1})} 
\leq C \| u_0 \|_{B^0_{\infty,1}}, 
\quad 
\| u_{n} \|_{\widetilde{L} ^1 (0,T; B^{1}_{\infty,1})} 
\leq 2\varepsilon, 
\\ \notag 
\| u_{n+1} - u_{n} \|
_{\widetilde L^2 (0,T ;  B^{-\delta}_{\infty,1})}
\leq C 2^{-n } .
\end{gather}
So we get a solution $u \in L^\infty (0,T ; B^0_{\infty,1} (\mathbb R)) 
\cap L^1 (0,T ; B^1_{\infty,1} (\mathbb R))$. 
Noting that $u_n$ is a bounded sequence in $L^\infty (\mathbb R)$ 
the dual of $L^1 (\mathbb R)$, we see that 
$u_n$ may converge to $u$ in the topology 
of dual weak sense in $L^\infty (\mathbb R)$ by taking a subsequence 
if necessary, so $u_n (t,x)$ converges to $u (t,x)$ for 
almost everywhere $x \in \mathbb R$. The Fatou lemma yields that 
$$
\| u(t) \|_{L^1} 
\leq \liminf _{n\to\infty} \| u_n(t) \|_{L^1} \leq \| u_0 \|_{L^1},
$$
where we have used the maximum principle in $L^1 (\mathbb R)$ 
assured by multiplying the equation by $u / |u|$, 
integrating and the integration by parts. 
This inequality proves $u(t) \in L^1 (\mathbb R)$, 
and we have that $u$ satisfies the integral equation \eqref{0723-6} 
thanks to the following estimate of nonlinear term
\begin{equation}\notag 
\begin{split}
\int_0^t \| e^{-(t-\tau)\Lambda} \partial_x u(\tau) ^2 
         \|_{L^1 \cap B^0_{\infty,1}}\, d\tau
\leq 
& 
C \int_0^t \| u \|_{L^1 \cap B^0_{\infty,1}} \| u \|_{B^1_{\infty,1}} \,d\tau
\\
\leq 
&C \| u \|_{L^\infty (0,T ; L^1 \cap B^0_{\infty,1})} 
   \| u \|_{L^1 (0,T ; B^1_{\infty,1})} , 
\end{split}
\end{equation}
although we omit  the detail. 
This estimates also implies the time continuity of $u$ in 
$L^1 (\mathbb R) \cap B^0_{\infty,1} (\mathbb R)$, 
since the linear part $e^{-t\Lambda} u_0$ is continuous 
in $L^1 (\mathbb R) \cap B^1_{\infty,1} (\mathbb R)$. 
The global regularity (see~\cite{KNS-2008,MW-2009}) and 
the smoothness by the analyticity yield that 
$u$ satisfies $u \in L^1_{loc} ([0,\infty) , \dot B^1_{\infty,1} (\mathbb R))$, 
which proves Lemma~\ref{lem:0723-1}. 
\end{pf}

\vskip3mm 

The following proposition is essential to handle 
the large time behavior of solutions. 

\begin{prop}\label{prop:global_bounds}
Let $u$ be a solution of the integral equation \eqref{0723-6} 
obtained in Lemma~\ref{lem:0723-1}. Then 
$u \in L^1 (0,\infty ; \dot B^1_{\infty,1} (\mathbb R))$. 

\end{prop}

\vskip3mm

\noindent 
{\bf Proof.} 
We start by proving the decay estimate of the solution in $L^\infty (\mathbb R)$ 
along the paper~\cite{CoCo-2004}. 
By taking $x_t$ such that $|u(t,x_t)| = \| u(t) \|_{L^\infty}$, 
we have from the equation that 
\begin{equation}\label{0807-1}
\partial_t \| u(t) \|_{L^\infty} 
\leq - \Lambda u(t,x_t) \operatorname{sgn} (u(t,x_t)). 
\end{equation}
Without loss of generality, we can assume that 
$u(t,x_t) \geq 0$, since it suffices to consider $-u(t,x_t)$ otherwise, 
so let $u(t,x_t) \geq 0$. We also put 
$$
B_\delta( x_t) 
:= \{ y \in \mathbb R \, | \, |y-x_t| \leq \delta \}, 
\quad 
\Omega _1 := 
\{ y \in B_\delta (x_t) \, | \, 
u(t,x_t) - u(t,y) \geq u(t,x_t ) /2 \},
$$
where $\delta$ will be taken later. Then it follows from 
$u(t,x_t) \geq u(t,y)$ for all $y \in \mathbb R$ that 
$$
\Lambda u(t,x_t) 
\geq 
P.V. \int_{\Omega _1} \frac{u(t,x_t) - u(t,y)}{|x_t -y|^2} \, dy 
\geq 
\frac{u(t,x_t)}{2 \delta ^2} \int_{\Omega _1} dy.
$$
On the other hand, 
we have from the maximum principle in $L^1 (\mathbb R)$ that 
$$
\| u_0 \|_{L^1} 
\geq \| u(t) \|_{L^1} 
\geq \int_{B_\delta (x_t) \setminus \Omega _1} |u(t,y) |dy
\geq \frac{u(t,x_t)}{2} 
   \Big(  \int_{B_\delta (x_t)} dy - \int_{\Omega _1} dy \Big) .
$$
The above inequalities, \eqref{0807-1} 
and $\int_{B_\delta (x_t)} dy = 2\delta$ imply that 
$$
\partial_t \| u(t) \| _{L^\infty}
\leq - \frac{u(t,x_t)}{2\delta^2} 
    \int_{\Omega_1} dy
\leq \frac{u(t,x_t)}{2\delta^2}
   \Big( \frac{2\| u_0 \|_{L^1}}{u(t,x_t)} - 2\delta \Big) .
$$
By taking $\delta = 2 \| u_0 \|_{L^1} /  u(t,x_t)$, we get that 
\begin{gather}\notag 
\partial_t \| u(t) \| _{L^\infty}
\leq - \frac{u(t,x_t)}{2\delta^2}
     \frac{ 2\| u_0 \|_{L^1}}{u(t,x_t)} 
= - \frac{\| u \|_{L^\infty} ^2}{4\| u_0 \|_{L^1}} ,
\\ \label{0807-2}
\| u(t) \|_{L^\infty} 
\leq \frac{\| u_0 \|_{L^\infty}}
   { 1 + \frac{\| u_0 \|_{L^\infty}}{4 \| u_0 \|_{L^1}} t} ,
\end{gather}
which completes the proof of time decay estimate of $u$ in $L^\infty (\mathbb R)$.

We turn to prove $u \in L^1 (0,\infty ; \dot B^1_{\infty,1} (\mathbb R))$, 
dividing into two cases; 
 $\| u_0 \|_{L^\infty} \leq \delta _0$, 
and $\| u_0 \|_{L^\infty} > \delta _0$, 
where $\delta _0$ is a small constant which will be taken 
later.

If $\| u_0 \|_{L^\infty} \leq \delta_0 $, 
the solution also satisfies $\| u(t) \|_{L^\infty} \leq \delta $ 
thanks to the maximum principle in $L^\infty (\mathbb R)$. 
Then we can estimate directly the integral equation 
by the maximal regularity estimate $\dot B^0_{\infty,1} (\mathbb R)$ 
and the bilinear estimate 
$\| u^2 \|_{\dot B^1_{\infty,1}} 
\leq \| u\|_{L^\infty} \| u \|_{\dot B^1_{\infty,1}}$ that 
\begin{equation}\notag 
\begin{split}
\| u \|_{L^\infty (0,\infty ; \dot B^0_{\infty,1}) 
  \cap L^1 (0,\infty ; \dot B^1_{\infty,1})} 
  \leq 
& 
C \| u_0 \|_{\dot B^0_{\infty,1}} 
   + C \int_0^\infty \| \partial _x u^2 \|_{\dot B^0_{\infty,1}} d\tau 
\\
\leq 
& C \| u_0 \|_{\dot B^0_{\infty,1}} 
   + C \|  u \|_{L^\infty (0,\infty; L^\infty)} 
      \| u \|_{L^1 (0,\infty ; \dot B^0_{\infty,1})}
\\
\leq 
& C \| u_0 \|_{\dot B^0_{\infty,1}} 
   + C \delta_0 
\| u \|_{L^\infty (0,\infty ; \dot B^0_{\infty,1}) 
  \cap L^1 (0,\infty ; \dot B^1_{\infty,1})} .
\end{split}
\end{equation}
By taking $\delta_0$ such that $C\delta _0\leq 1/2$, 
we have that 
$$
\| u \|_{L^\infty (0,\infty ; \dot B^0_{\infty,1}) 
  \cap L^1 (0,\infty ; \dot B^1_{\infty,1})} 
  \leq 
2 C \| u_0 \|_{\dot B^0_{\infty,1}} ,
$$
which proves 
$u \in L^1(0,\infty ; \dot B^1_{\infty,1} (\mathbb R))$.

If $\| u_0 \|_{L^\infty} > \delta _0$, we apply the estimate \eqref{0807-2} 
to get that 
$$
\| u(t_0) \|_{L^\infty} \leq \delta _0 
\quad \text{for } 
t_0 = \frac{4\| u_0 \|_{L^1}}{\| u_0 \|_{L^\infty}} 
 \Big( \frac{\| u_0 \|_{L^\infty}}{\delta_0} -1 \Big) . 
$$
Hence, the previous case of small initial data implies that 
$u \in L^1 (t_0 , \infty ; \dot B^1_{\infty,1} (\mathbb R))$, 
which assures $u \in L^1 (0,\infty ; \dot B^1_{\infty,1} (\mathbb R))$. 
\hfill 
$\Box$

\vskip3mm 

We next prove the decay estimates of solutions. 

\begin{prop}\label{prop:decay}
Let $u$ be a solution of the integral equation \eqref{0723-6} 
obtained in Lemma~\ref{lem:0723-1} 
such that 
$u \in C([0,\infty), \dot B^s_{\infty,1}(\mathbb R))$ 
for any $s > 0$. 
Let $\alpha > 0$ and $1 \leq p \leq \infty$. 
Then
\begin{gather}\label{0723-10}
\| u(t) \|_{L^p} \leq C t^{-(1-\frac{1}{p})} 
\quad \text{for any } t \geq 1.
\\ \label{0723-11}
\int_0^\infty \| t^\alpha u(t) \|_{\dot B^{1+\alpha}_{\infty,1}} 
\, dt < \infty
\quad \text{and} \quad 
\| |\nabla|^\alpha u(t) \|_{L^p} 
\leq C t^{-\alpha - (1-\frac{1}{p})} 
\quad \text{for any } t \geq 1.
\end{gather}

\end{prop}

\begin{pf}
The decay estimate \eqref{0723-10} is a consequence of 
\eqref{0807-2} or 
Proposition~5.2 in \cite{Iw-2015}, 
since $u \in L^1 (0,\infty; \dot B^1_{\infty,} (\mathbb R))$ 
by Proposition~\ref{prop:global_bounds}. 
We also notice that the decay estimate of 
$\| |\nabla|^\alpha u(t) \|_{L^\infty}$ in \eqref{0723-11} 
can be proved by Proposition~5.2 in \cite{Iw-2015} 
once the integrability in \eqref{0723-11} is obtained. 
Hence, all we have to do is to prove the former part of \eqref{0723-11}. 

We prove the integrability in \eqref{0723-11}.
Since $u$ is smooth, it is sufficient to consider the 
integrabiity for large $t$. By the decay estimate \eqref{0723-10}, 
we see that for any $\delta > 0$, there exists $t_\delta > 0$ 
such that $\| u(t) \|_{L^\infty} \leq \delta $ for any 
$t \geq t_\delta$, so we may consider the integrability 
on the time interval with smallness of $L^\infty$ norm. 
Put 
$$
v(t) := u(t_\delta + t) \quad \text{for } t \geq 0.
$$

For any $T > 0$, we estimate the integral equation for $v$ 
similarly to the proof of Proposition~5.2 in \cite{Iw-2015}. 
For the linear part, it follows from the smoothing effect 
and the maximal regularity for $e^{-t\Lambda}$ that 
$$
\int_0^T \| t^\alpha e^{-t\Lambda} v(0) 
\|_{\dot B^{1+\alpha}_{\infty,1}} 
\, dt
\leq 
C \int_0^T \|  e^{-\frac{t}{2}\Lambda} v(0) 
\|_{\dot B^{1}_{\infty,1}} 
\, dt 
\leq C \| v(0) \|_{\dot B^0_{\infty,1}} .
$$
As to the nonlinear part, 
we decompose $[0,t]$ into $[0,t/2]$ and $[t/2,t]$ to have that 
\begin{equation}\notag 
\begin{split}
& \int_0^T
\Big\| t^\alpha 
\int_0^{\frac{t}{2}} e^{-(t-\tau)\Lambda} \partial_x v(\tau)^2 
\, d\tau
\Big\|_{\dot B^{1+\alpha}_{\infty,1}} dt
\\
\leq 
& 
C \int_0^T t^\alpha 
  \int_0^{\frac{t}{2}} (t-\tau)^{-\alpha}
    \| e^{-\frac{t-\tau}{2} \Lambda} v(\tau)^2 \|_{\dot B^2_{\infty,1}} 
    d\tau \, dt
\leq 
C \int_0^T
  \int_{2 \tau}  ^T
    \| e^{-\frac{t-\tau}{2} \Lambda} v(\tau)^2 \|_{\dot B^2_{\infty,1}} 
 dt \, d\tau 
\\
\leq 
& 
C \int_0^T 
    \| v(\tau)^2 \|_{\dot B^1_{\infty,1}} 
\, d\tau 
\leq 
C \| v \|_{L^\infty (0,T; L^\infty)} 
  \| v \|_{L^1 (0,T; \dot B^1_{\infty,1})},
\end{split}
\end{equation}
where we have used the smoothing effect and the maximal regularity 
estimate for $e^{-\frac{t-\tau}{2}\Lambda}$ 
and the bilinear estimate 
$\| v^2 \|_{\dot B^1_{\infty,1}} 
\leq \| v \|_{L^\infty} \| v \|_{\dot B^1_{\infty1}}$. 
On the integral on $[t/2 , t]$, 
it follows from $t^\alpha \leq 2^\alpha \tau ^{\alpha}$ for 
$t/2 \leq \tau \leq t$, the maximal regularity estimate 
in $\dot B^{\alpha+1}_{\infty,1} (\mathbb R)$ 
and the bilinear estimate 
$\| v^2 \|_{\dot B^{1+\alpha}_{\infty,1}} 
\leq C\| v \|_{L^\infty} \| v \|_{\dot B^{1+\alpha}_{\infty,1}}$ 
that 
\begin{equation}\notag 
\begin{split}
\int_0^T
\Big\| t^\alpha 
\int_{\frac{t}{2}}^t e^{-(t-\tau)\Lambda} \partial_x v(\tau)^2 
\, d\tau
\Big\|_{\dot B^{1+\alpha}_{\infty,1}} dt
\leq 
& 
C 
\int_0^T
\int_{\frac{t}{2}}^t 
\| e^{-(t-\tau)\Lambda} \tau ^{\alpha}  v(\tau)^2 
\|_{\dot B^{2+\alpha}_{\infty,1}}
\, d\tau \,dt
\\
\leq 
& 
C 
\int_0^T
\int_{\tau}^T 
\| e^{-(t-\tau)\Lambda} \tau ^{\alpha}  v(\tau)^2 
\|_{\dot B^{2+\alpha}_{\infty,1}}
\, dt \,d\tau
\\
\leq 
& 
C 
\int_0^T
\| \tau ^{\alpha}  v(\tau)^2 
\|_{\dot B^{1+\alpha}_{\infty,1}}
\, dt \,d\tau
\\
\leq 
& 
C 
\| v \|_{L^\infty (0,T; L^\infty)} 
\int_0^T
\tau ^{\alpha} \| v \|_{\dot B^{1+\alpha}_{\infty,1}}
\,d\tau .
\end{split}
\end{equation}
We also have on the norm of 
$L^1 (0,T ; \dot B^1_{\infty,1} (\mathbb R))$ 
from the maximal regularity estimate in 
$\dot B^0_{\infty,1} (\mathbb R)$ and the bilinear estimate 
$\| v^2 \|_{\dot B^1_{\infty,1}} 
\leq C\| v \|_{L^\infty} \| v \|_{\dot B^1_{\infty,1}}$ 
that 
$$
\| v \|_{L^1 (0,T ; \dot B^1_{\infty,1})} 
\leq C \| v(0) \|_{\dot B^0_{\infty,1}}
+ C \| v \|_{L^\infty (0,T ; L^\infty)} 
    \| v \|_{L^1 (0,T ; \dot B^1_{\infty,1})} . 
$$
The above four estimates and 
$\| v(t) \|_{L^\infty} \leq \delta$ yield that 
\begin{equation}
\begin{split}
& 
\| t^\alpha v(t) \|_{L^1_t (0,T ; \dot B^{1+\alpha}_{\infty,1})}
+ \| u(t) \|_{L^1(0,T ;\dot B^1_{\infty,1})} 
\\
\leq 
& C \| v(0) \|_{\dot B^0_{\infty,1}}
+ C \delta 
\big( \| t^\alpha v(t) \|_{L^1_t (0,T ; \dot B^{1+\alpha}_{\infty,1})}
    + \| u(t) \|_{L^1(0,T ;\dot B^1_{\infty,1})} 
\big) .
\end{split}
\end{equation}
Here taking $\delta $ such that $C \delta \leq 1/2$, 
where $C$ is a constant appearing in the above estimate, 
we obtain 
$$
\| t^\alpha v(t) \|_{L^1_t (0,T ; \dot B^{1+\alpha}_{\infty,1})}
+ \| u(t) \|_{L^1(0,T ;\dot B^1_{\infty,1})} 
\leq 2 C \|  v(0) \|_{\dot B^0_{\infty,1}} 
\quad \text{for any } T>0.
$$
Hence the integrability in \eqref{0723-11} is verified 
and we finish to prove Proposition \ref{prop:decay}. 
\end{pf}

\vskip3mm

Based on the lemma and the propositions, 
we prove the large time behavior. 

\vskip3mm 

\noindent 
{\bf Proof of large time behavior \eqref{largetime}}. 
Let $u$ be a global solution, which is obtained 
by {\rm (i)} of Theorem~\ref{thm:Burgers} and the global regularity, 
with initial data $u_0$ 
satisfying \eqref{B_initial} and $u_0 \in L^1 (\mathbb R)$. 
We see 
that $u(t) \in B^s_{\infty,1} (\mathbb R)$ for any $s \geq 0$, $t > 0$, 
since we can have regularity 
as much as it is needed thanks to the analyticity and the global regularity.  
Once $u(t) \in L^1 (\mathbb R)$ for all $t$ near $0$ is proved, 
we are able to apply Lemma~\ref{lem:0723-1} and 
Propositions~\ref{prop:decay} by regarding the initial data 
as $u(t_0)$ for time $t_0 > 0$ near $0$ 
to obtain the decay estimates 
\eqref{0723-10} and \eqref{0723-11}, 
which proves the large time behavior \eqref{largetime} 
in the same argument of the paper~\cite{Iw-2015} 
(see the proof of (1.4)). 
All we have to do is to prove that 
$u(t) \in L^1 (\mathbb R)$ for all $t$ near $0$.

Put $u_{0,N} := S_N u_0$, where 
$S_N = (\psi + \sum _{j=1} ^N \phi_j )*$. 
We denote by $u_N$, $u$ the solutions with the initial data 
$u_{0,N}$, $u_{0}$, respectively. 
It follows from 
$u_{0,N} \in L^1 (\mathbb R) \cap B^0_{\infty,1} (\mathbb R)$ 
and Lemma~\ref{lem:0723-1} that $u_N$ satisfies the energy 
identity in $L^2 (\mathbb R)$
$$
\| u_N(t) \|_{L^2}^2 
+ 2 \int_0^t \| \Lambda^{\frac{1}{2}} u_N \|_{L^2} ^2 d\tau
\leq \| u_{0,N} \| _{L^2}^2. 
$$
and hence, 
$$
\| u_N(t) \|_{L^2}^2 
\leq \| u_{0,N} \| _{L^2}^2 
\leq \| u_0 \|_{L^2} ^2 .
$$
Since $L^2 (\mathbb R)$ is a Hilbert space, 
we can take a subsequence of $\{ u_N(t) \}_{N=1}^\infty$, 
denoted by the same, $u_N (t)$ converges to an element 
$\widetilde u(t)$ of 
$L^2(\mathbb R)$ in the weak topology of $L^2 (\mathbb R)$. 
On the other hand, we can prove that $u_N(t)$ tends to 
$u(t)$ in $\mathcal S'(\mathbb R)$ for almost every $t$. 
In fact, it follows from the same proof of \eqref{0402-d} that 
for $q = \infty$
\begin{equation}\label{0723-20}
\begin{split}
& 
\| u_{N} - u \|
_{\widetilde L^2 (0,T ;  B^{-\delta}_{\infty,q})}
\\
\leq 
& 
C \| e^{-tc \Lambda} (u_{0,N} - u_0) 
  \|_{\widetilde L^2 (0,T ;  B^{-\delta}_{\infty,q})} 
\\
&
+ C 
 \big( \| u_N \|_{\widetilde L^2 (0,T ;  B^{\frac{1}{2}}_{\infty,\infty})}
      +\| u \|_{\widetilde L^2 (0,T ;  B^{\frac{1}{2}}_{\infty,\infty})}
 \big)     
  \| u_{N} - u 
  \|_{\widetilde L^2 (0,T ;  B^{-\delta}_{\infty,q} )} .
\end{split}
\end{equation}
We note that 
$$
\lim_{T \to 0} 
\Big(
  \sup _{N} 
 \| u_N \|_{\widetilde L^2 (0,T ;  B^{\frac{1}{2}}_{\infty,\infty})}
      +\| u \|_{\widetilde L^2 (0,T ;  B^{\frac{1}{2}}_{\infty,\infty})}
\Big) 
= 0
$$
which is assured by that $u_{0,N}$ is defined by restricting 
the frequency of $u_0$. 
Hence, by taking $T= T_0 > 0$ sufficiently small, the inequality \eqref{0723-20} 
yields that 
$$
\| u_{N} - u \|
_{\widetilde L^2 (0,T_0 ;  B^{-\delta}_{\infty,\infty})}
\leq 2C \| e^{-tc \Lambda} (S_{N}u_0 - u_0) 
  \|_{\widetilde L^2 (0,T_0 ;  B^{-\delta}_{\infty,\infty})} 
\leq 2C T_0^{\frac{1}{2}} \| u_{0,N} - u_0 \|_{B^{-\delta}_{\infty,\infty}}
\to 0
$$
as $N \to \infty$. 
So $u_N(t,x)$ converges to $u(t,x)$ as $N \to \infty$ 
in $\mathcal D' ((0,T_0)\times \mathbb R )$. 
Therefore for almost every $t$, 
$u(t) = \widetilde u(t)$ in $\mathcal S '(\mathbb R)$ 
and $u_N(t,x)$ tends to $u(t,x)$ for almost every $x \in \mathbb R$ 
by the weak convergence in $L^2 (\mathbb R)$ of $u_N(t)$ to $u(t)$. 
By applying the Fatou Lemma, we have that 
$$
\| u(t) \|_{L^1 } 
\leq \liminf_{N \to \infty} \| u_N (t) \|_{L^1}
\leq \liminf_{N \to \infty} \| u_{0,N} \|_{L^1}
\leq C \| u_0 \|_{L^1} < \infty, 
$$
where we have used the maximum principle in $L^1 (\mathbb R)$, 
namely, $\| u_{N} (t) \|_{L^1} \leq \| u_{0,N} \|_{L^1}$ 
obtained by multiplying the equation by $u/|u|$, integrating 
and the integration by parts. 
Therefore, $u(t) \in L^1 (\mathbb R)$ for almost every $t \leq T_0$ 
and the smoothness 
of $u$ assures that $u(t) \in L^1 (\mathbb R)$ for all $t \leq T_0$. 
We complete the proof of the large time behavior 
in Theorem~\ref{thm:Burgers}. 
\hfill $\Box$

\vskip3mm 

\noindent
{\bf Acknowledgements.}
The author was supported by the Grant-in-Aid for Young Scientists (A) (No.~17H04824)
from JSPS 
and by Advancing Strategic International Networks to Accelerate 
the Circulation of Talented Researchers from JSPS.

\vskip3mm 
%
%
%
%

\end{document}